\documentclass[a4paper, reqno]{amsart}
\usepackage{enumerate}
\usepackage{amsthm}
\usepackage{amsmath}
\usepackage[mathscr]{euscript}
\usepackage{amssymb}
\usepackage[noadjust]{cite}
\usepackage{comment}
\usepackage{color}
\usepackage{filecontents}
\usepackage{enumerate}

\usepackage[plainpages=false,colorlinks,hyperindex,pdfpagemode=None,bookmarksopen,linkcolor=blue,citecolor=blue,urlcolor=blue]{hyperref}

\newtheorem{thm2}{Theorem}
\newtheorem{prop2}[thm2]{Proposition}
\newtheorem{prob2}[thm2]{Problem}

\newtheorem{thm}{Theorem}[section]
\newtheorem{lem}[thm]{Lemma}
\newtheorem{prop}[thm]{Proposition}
\newtheorem{cor}[thm]{Corollary}

\theoremstyle{definition}

\newtheorem{defn}[thm]{Definition}

\theoremstyle{remark}
\newtheorem{eg}[thm]{Example}
\newtheorem{rem}[thm]{Remark}

\numberwithin{equation}{section}

\newcommand{\ti}{\tilde}

\newcommand{\smnoind}{\smallskip\noindent}

\newcommand{\ps}{\mathcal{PS}}
\newcommand{\id}{\mathrm{id}}

\newcommand{\BC}{\mathbb{C}}
\newcommand{\BN}{\mathbb{N}}
\newcommand{\BM}{\mathbb{M}}

\newcommand{\CO}{\mathcal{O}}

\newcommand{\CL}{\mathcal{L}}

\newcommand{\CS}{\mathcal{S}}
\newcommand{\CU}{\mathcal{U}}

\newcommand{\CT}{\mathcal{T}}
\newcommand{\D}{\mathcal{D}}

\newcommand{\CH}{\mathcal{H}}
\newcommand{\CK}{\mathcal{K}}

\newcommand{\Z}{\mathrm{Z}}

\newcommand{\bs}{\mathbf{s}}
\newcommand{\bp}{\mathbf{p}}
\newcommand{\ks}{\mathfrak{s}}
\newcommand{\bz}{\mathbf{z}}
\newcommand{\sa}{\mathrm{sa}}

\newcommand{\ba}{\mathbf{a}}

\newcommand{\BI}{\mathbf{I}}

\newcommand{\KI}{\mathfrak{I}}

\newcommand{\m}{{\min}}

\newcommand{\I}{\mathrm{I}}

\newcommand{\RB}{\mathbf{B}}
\newcommand{\RS}{\mathbf{S}}
\newcommand{\CC}{\mathcal{C}}

\newcommand{\BP}{\mathbf{P}}

\newcommand{\de}{\mathbf{1}}

\newcommand{\cdd}{ }

\newcommand{\ida}[1]{{1_{#1^{**}_\ba}}}

\begin{document}

\title[Reconstructing the Jordan structure]{Reconstructing the Jordan structure of $C^*$-algebras from the metric on positive unit spheres}

	\author[Leung, Ng and Wong]{Chi-Wai Leung
		\and Chi-Keung Ng \and Ngai-Ching Wong}
	
	\address[Chi-Wai Leung]{Department of Mathematics, The Chinese University of Hong Kong, Hong Kong, China}
	\email{cwleung@math.cuhk.edu.hk}
	
	\address[Chi-Keung Ng]{Chern Institute of Mathematics and LPMC, Nankai University, Tianjin 300071, China.}
	\email{ckng@nankai.edu.cn; ckngmath@hotmail.com}

	\address[Ngai-Ching Wong]{Department of Applied Mathematics, National Sun Yat-sen University, Kaohsiung, 80424, Taiwan.}
	\email{wong@math.nsysu.edu.tw}

\thanks{Corresponding author: Ngai-Ching Wong, wong\@math.nsysu.edu.tw}



\date{\today}

\begin{abstract}
Let $A$ and $B$ be unital $C^*$-algebras, and let
$\RS_{A^+}$ and $\RS_{B^+}$ denote their positive unit spheres.
We first show that the diametral relation on the positive unit sphere
already determines Jordan structure.
More precisely, suppose that
$\Lambda:\RS_{A^+}\to\RS_{B^+}$ is a bijection satisfying
\[
   \|\Lambda(u)-\Lambda(v)\|=1
   \quad\Longleftrightarrow\quad
   \|u-v\|=1
   \qquad (u,v\in\RS_{A^+}).
\]
If $A$ has no two-dimensional irreducible $^*$-representation, then
$A$ and $B$ are Jordan $^*$-isomorphic. Moreover, the Jordan
$^*$-isomorphism can be chosen to be compatible with $\Lambda$ on the
peak and zero projections of positive norm-one elements. The proof is
based on recovering the order and orthogonality structure of Akemann
closed projections in the bidual.

As a main application, we give an affirmative answer to the positive-sphere extension problem for unital $C^*$-algebras of real rank
zero.
Namely, if $A$ has real rank zero (but could has two-dimensional irreducible $^*$-representations) and
$\Lambda:\RS_{A^+}\to\RS_{B^+}$ is a surjective isometry, then
$\Lambda$ extends to a Jordan $^*$-isomorphism from $A$ onto $B$.
\end{abstract}

\maketitle

\section{Introduction}

In this paper, we continue the study of how much geometric information
is needed to recover the Jordan $^*$-algebra structure of a
$C^*$-algebra; see, for example,
\cite{Tana,FP,FP2,FP3,Mori,Peralta18,LN,LNW,MO,Kad51,leung2021variant}.
One of our conclusions is that surprisingly weak metric information on the
positive unit sphere already encodes substantial algebraic structure.

More precisely, we show that for a unital $C^*$-algebra with no
two-dimensional irreducible $^*$-representation, the Jordan
$^*$-algebra structure can be recovered from the single binary
relation of being at distance one, on the positive unit sphere.
We also prove that every surjective
isometry between the positive unit spheres of unital $C^*$-algebras
extends to a  Jordan $^*$-isomorphism whenever the domain algebra has
real rank zero. These results continue the programme initiated in
\cite{LNW} of recovering $C^*$-algebras, up to Jordan
$^*$-isomorphism, from restricted metric data.

For a $C^*$-algebra $A$, let
\[
   \RS_{A^+}=\{a\in A:a\geq 0,\ \|a\|=1\}
\]
denote its \emph{positive unit sphere}. Since its elements are positive
contractions, the distance between any two elements of $\RS_{A^+}$ is
at most one. We therefore consider the \emph{diametral relation}
\begin{equation}\label{eqt:def-diam-rel}
   \D_{\RS_{A^+}}
   =
   \{(a,b)\in\RS_{A^+}\times\RS_{A^+}:\|a-b\|=1\}.
\end{equation}
Preserving this relation is substantially weaker than preserving the
full metric. Nevertheless, this information is sufficient to recover
the Jordan structure under a mild representation-theoretic assumption, as we see in the following result (which will be established in Theorem \ref{thm:pre-diam-rel}).
Throughout this paper, all Jordan $^*$-isomorphisms are assumed to be complex linear.

\begin{thm2}\label{thm:main1}
Let $A$ and $B$ be unital $C^*$-algebras.
Suppose that $A$ has no
two-dimensional irreducible $^*$-representation. If a bijection
$\Lambda:\RS_{A^+}\to\RS_{B^+}$ satisfies
\begin{equation}\label{eqt:pre-diam}
   \|a-b\|=1
   \quad\Longleftrightarrow\quad
   \|\Lambda(a)-\Lambda(b)\|=1
   \qquad(a,b\in\RS_{A^+}),
\end{equation}
then there is a Jordan $^*$-isomorphism $\Psi:A\to B$ such that,
for each $a\in\RS_{A^+}$,
\[
   \Psi^{**}(\bp(a))=\bp(\Lambda(a)),
   \qquad
   \Psi^{**}(\bz(a))=\bz(\Lambda(a)),
\]
where
\[
   \bp(a)=w^*\text{-}\lim_n a^n,
   \qquad
   \bz(a)=w^*\text{-}\lim_n(1-a)^n
\]
are the spectral projections in $A^{**}$ corresponding, respectively,
to the peak and zero parts of $a$. In particular,
\[
   \Psi|_{\BP_A\setminus\{0\}}
   =
   \Lambda|_{\BP_A\setminus\{0\}}.
\]
\end{thm2}

Theorem \ref{thm:main1} shows that the
diametral relation determines the Jordan structure of the algebra.
It does not, however, assert that the given bijection $\Lambda$ itself
extends linearly to the whole algebra. Indeed, a bijection preserving
the diametral relation need not admit any real linear extension; see
Example \ref{eg:counter-eg-lin-ext} and also
\cite[Example 3.3]{LNW}.
Thus, the geometry encoded by the relation
$\|a-b\|=1$ can determine the Jordan structure even when the map
preserving that relation is not itself linearizable.

Theorem \ref{thm:main1} gives an affirmative answer to the following
question from \cite[Problem 4.10]{LNW} for unital $C^*$-algebras
without two-dimensional irreducible $^*$-representa\-tions. The problem
remains open in the presence of such representations.

\begin{prob2}\label{prob:weak-Tingley}
Let $A$ and $B$ be $C^*$-algebras. If there exists a bijection
$\Lambda:\RS_{A^+}\to\RS_{B^+}$ preserving the diametral relation in
both directions, must $A$ and $B$ be Jordan $^*$-isomorphic?
\end{prob2}

The preceding theorem recovers the Jordan structure from the
diametral relation, but does not linearize the given map.
Our second
main result as stated below (see Theorem \ref{thm:main}) addresses precisely this stronger extension problem.
In contrast with Theorem \ref{thm:main1}, this theorem requires no restriction on
two-dimensional irreducible representations, but it uses Theorem \ref{thm:main1} as a main step in its proof.

\begin{thm2}\label{thm:main2}
Let $A$ and $B$ be unital $C^*$-algebras.
Suppose that $A$ has
real rank zero. Every surjective isometry
\[
\Lambda:\RS_{A^+}\to\RS_{B^+}
\]
extends to a Jordan $^*$-isomorphism from $A$ onto $B$.
\end{thm2}

Notice that only the domain algebra is assumed to have real rank zero.
Theorem \ref{thm:main2} gives an affirmative answer, in the real rank zero
setting, to the following \emph{positive-sphere extension problem},
formulated in \cite{Peralta18,LN}.

\begin{prob2}\label{prob:Tingley}
Let $A$ and $B$ be $C^*$-algebras. Does every surjective isometry
\[
\Lambda:\RS_{A^+}\to\RS_{B^+}
\]
extend to a  Jordan $^*$-isomorphism from $A$ onto $B$?
\end{prob2}

Affirmative answers to Problem \ref{prob:Tingley} were previously known
for algebras of bounded operators on Hilbert spaces
(\cite{Na13,nagy2018isometries,peralta2019unit}), commutative
$C^*$-algebras (\cite{leung2021variant,ELM26}), finite direct sums of
homogeneous type $\BI$ $W^*$-algebras (\cite{LN}), and general
$AW^*$-algebras (\cite{LNW}). A corresponding result for
$JBW^*$-algebras whose type $\BI_2$ parts are zero or atomic has also
been obtained in \cite{PS26}.

Problem \ref{prob:Tingley} belongs to the broader study of the
algebraic information encoded by norm geometry.
In fact, the Mazur--Ulam
theorem \cite{MU} asserts that every surjective isometry between real
normed spaces is affine. Mankiewicz \cite{mankiewicz1972extension}
extended this principle to isometries between open connected subsets
and, more generally, between convex bodies.

Alone this line, Tingley's problem \cite{Tingley} asks whether every surjective
isometry between the unit spheres of real Banach spaces extends to a
real linear isometry of the ambient spaces; see, for example,
\cite{CD,KM,Peralta18,BJM} for developments and related results. In
the setting of $C^*$-algebras, following work of Tanaka, Peralta,
Mori, and others (\cite{Tana,FP,FP2,FP3,Mori,Peralta18}), Mori and
Ozawa proved  in \cite{MO} that every surjective isometry from the unit
sphere of a $C^*$-algebra onto the unit sphere of a normed space
extends to a real linear isometry.
Very recently, a preprint by OpenAI
announced an affirmative solution of Tingley's problem for arbitrary
real Banach spaces \cite{OpenAI-Tingley}.

Between $C^*$-algebras, an isometric linear extension as above is a real $JB^*$-triple isomorphism
(\cite[Theorem 6.4]{CDRV}).
On the other hand, Kadison's isometry theorem (\cite{Kad51}) states
that every surjective complex linear isometry between unital
$C^*$-algebras has the form
\[
a\mapsto uJ(a),
\]
where $u$ is a unitary and $J$ is a Jordan $^*$-isomorphism.

When the
objective is to recover the complex Jordan $^*$-algebra structure,
the positive unit sphere is a natural object to consider.
Indeed, any real linear isometric extension of a surjective isometry
between positive unit spheres is an order isomorphism on the
self-adjoint parts, and Kadison's order-theoretic characterization
(\cite{Kad52}) then yields a Jordan isomorphism whose complexification
is a Jordan $^*$-isomorphism. Conversely, every Jordan
$^*$-isomorphism restricts to a surjective isometry between positive
unit spheres. This is the motivation behind Problem
\ref{prob:Tingley}.

For noncommutative $C^*$-algebras, the most general result available
prior to the present work was the $AW^*$-algebra case treated in
\cite{LNW}. The principal difficulty in passing from $AW^*$-algebras
to general $C^*$-algebras is the possible scarcity of projections in
the algebra itself. In the $AW^*$-setting, the argument in \cite{LNW}
recovers the ortho-lattice of projections from the diametral relation.
Hamhalter's extension (\cite{Ham}) of Dye's theorem
(\cite{dye1955geometry}) then recovers the Jordan structure and gives a
positive answer to Problem \ref{prob:weak-Tingley} for $AW^*$-algebras
without a type $\BI_2$ summand.
This is the first step to the solution to Problem
\ref{prob:Tingley} in the $AW^*$-case.

A general $C^*$-algebra, however, may contain no nontrivial
projections, so this approach cannot be applied directly. We overcome
this difficulty by using Akemann's noncommutative topology
(\cite{Ake69,Ake71,ORT}). The appropriate substitute for the projection
lattice is the set $\CC(A)$ of closed projections in $A^{**}$ relative
to $A$. Although these projections need not belong to $A$ itself,
their order and orthogonality can still be recovered from the
diametral relation on $\RS_{A^+}$. More precisely, we prove the
following result (see Proposition \ref{prop:pre-diam-rel}).

\begin{prop2}\label{prop2:bi-orth-pre}
Let $A$ and $B$ be unital $C^*$-algebras, and let
$\Lambda:\RS_{A^+}\to\RS_{B^+}$ be a bijection satisfying
\eqref{eqt:pre-diam}. Then there exists a bi-orthogonality preserving
order isomorphism
\[
\Gamma:\CC(A)\to\CC(B)
\]
such that
\[
   \Gamma(\bp(a))=\bp(\Lambda(a)),
   \qquad
   \Gamma(\bz(a))=\bz(\Lambda(a))
   \qquad(a\in\RS_{A^+}).
\]
\end{prop2}

A key ingredient in the construction of $\Gamma$ is a noncommutative Urysohn
lemma for separable $C^*$-algebras, deduced from
\cite[Theorem 1.1]{Ake}; see Corollary \ref{cor:Ury}. Once the closed
projections and their orthogonality structure have been recovered, an
analogue of Dye's theorem, first obtained in \cite{DN}, provides the
passage from this projection-theoretic information to the Jordan
structure.

\begin{prop2}[Dye's theorem for closed projections]
\label{prop:Dye-closed-proj}
Let $A$ and $B$ be unital $C^*$-algebras, and suppose that $A$ has no
two-dimensional irreducible $^*$-representation. For every
bi-orthogonality preserving bijection
\[
\Gamma:\CC(A)\to\CC(B),
\]
there exists a Jordan $^*$-isomorphism $\Psi:A\to B$ such that
\[
   \Gamma=\Psi^{**}|_{\CC(A)}.
\]
\end{prop2}

For completeness, Proposition \ref{prop:DN-3.8} gives a self-contained
proof of this result without using the language of ortho-sets.
Combining it with Proposition \ref{prop2:bi-orth-pre} yields
Theorem \ref{thm:main1}.

Theorem \ref{thm:main2} addresses the additional task of extending
the given map $\Lambda$ itself. Here the full metric structure of the
positive unit sphere is used, rather than merely the diametral relation.
The real rank zero assumption supplies finite-spectrum approximations helping to pass from the metric information
encoded by $\Lambda$ to a global Jordan $^*$-isomorphism.
Importantly,
this tactic can also be used to deal with algebras with two-dimensional irreducible
$^*$-representations (although a separated argument is needed for this case).

The paper is organized as follows. Section \ref{S:projections}
develops the algebraic and order-theoretic information encoded by the
diametral relation, including Proposition \ref{prop:order-1} and
Corollaries \ref{cor:v-in-hered} and \ref{cor:met-char-proj}.
Section \ref{S:weakT} provides the important result concerning connection of  closed projections (Proposition \ref{prop2:bi-orth-pre}), and proves
the Jordan structure recovering theorem
(Theorem \ref{thm:main1}).
Section \ref{S:T} verifies the
extension theorem for positive-sphere isometries when the domain
algebra has real rank zero (Theorem \ref{thm:main2}). Finally,
Appendix \ref{sec:open-proj} collects the required preliminary results
on open and closed projections.

\section{Notations and Preliminaries}\label{S:prel}

Throughout this paper, $A$ is a (complex) unital $C^*$-algebra with   identity  $1_A$, or simply  $1$.
We write $A^*$ and $A^{**}$ for the Banach dual space and bidual space of $A$, respectively.
Denote by $A_\sa$ and $A^+$ the self-adjoint part and the positive part of $A$, and by $\RB_A$ and $\RS_A$ the closed unit ball and
the unit sphere of $A$, respectively.
The \emph{positive closed unit ball} and the \emph{positive unit sphere} of $A$ are, respectively, the  metric spaces (under the norm-metric)
$$
\RB_{A^+}:= \{a\in A^+: \|a\|\leq 1 \} \quad\text{and}\quad
\RS_{A^+}:= \{a\in A^+: \|a\|= 1 \}.
$$
Let
$\Z(A)$ be the center of $A$, and
 $\BP_A$ be the set of all projections in $A$, respectively.
We use $A^{-1}$ to denote the   set of invertible elements of $A$, and put
$$
\RB_{A^+}^{-1} := \RB_{A^+}\cap A^{-1} \quad\text{and}\quad \RS_{A^+}^{-1} := \RS_{A^+}\cap A^{-1}.
$$
We also use $\ps(A)$ to denote  the set of all pure states on $A$.

\smnoind

We consider the canonical von Neumann algebra structure on $A^{**}$ and
regard
\begin{align*}
	A\subseteq A^{**}
\end{align*}
as a $C^*$-subalgebra.
For a subset $E\subseteq \BP_{A^{**}}$ and $p,q\in \BP_{A^{**}}$, the notations
$$\bigvee E\quad \text{and} \quad p\vee q$$
stand for the suprema taken in the complete lattice $\BP_{A^{**}}$, while the notations
$$
\bigwedge E\quad \text{and} \quad p\wedge q
$$
stand for the infima taken in $\BP_{A^{**}}$.

\smnoind

Let $B$ be another unital $C^*$-algebra.
If $\CS\subseteq \BP_{A^{**}}$ and $\CT\subseteq \BP_{B^{**}}$, then a bijection $\Gamma:\CS\to \CT$ is said to be \emph{bi-orthogonality preserving} if
for every $r,s\in \CS$, one has
\begin{equation}\label{eqt:def-bi-ortho-pres}
\Gamma(r)\Gamma(s) = 0 \quad \text{if and only if}\quad rs=0.
\end{equation}

\smnoind

The following lemma lists some known facts  (see, e.g., \cite{LNW}).

\begin{lem}\label{lem:peak-and-zero}
	Let $A$ be a unital $C^*$-algebra,
	$p\in \BP_A\setminus \{0\}$, $r\in \BP_{A^{**}}$ and $a,b\in \RS_{A^+}$.
	
	\smnoind
	(a) $\|a-b\| = 1$ if and only if there exists $\omega\in \ps(A)$ such that $\{\omega(a), \omega(b)\} = \{0,1\}$.
	
	\smnoind
	(b) $a\in A^{-1}$ if and only if $\omega(a)\neq 0$ for all $\omega\in\ps(A)$.
	
	\smnoind
	(c) If $\omega\in \ps(A)$ with $\omega(r)=1$, then $\omega(rar) = \omega(ar) =\omega(ra) =\omega(a)$.

	\smnoind
	(d) $\big\{c\in \RS_{A^+}^{-1}: p\leq c\big\} = p +  \RB_{((1-p)A(1-p))^+}^{-1}$.
\end{lem}

Let $E\subseteq \RS_{A^+}$.
We put
\begin{equation}\label{eqt:def-orth-comp}
	E^\de := \{u\in \RS_{A^+}: \|u-e\|=1,  \text{ for every }e\in E \}.
\end{equation}
It is obvious that $E\subseteq (E^\de)^\de$, and we have
\begin{equation}\label{eqt:111}
	((E^\de)^\de)^\de  = E^\de.
\end{equation}

\begin{lem}[{\cite[Lemma 2.2]{LNW}}]\label{lem:sum-peak}
	Let $A$ be a unital $C^*$-algebra.

	\smnoind
	(a) A projection $p\in \BP_A\setminus \{0\}$ is central if and only if $p\in (\BP_A\setminus \{p\})^\de$. 	
	
	\smnoind
	(b) If $p\in \BP_{\Z(A)}\setminus \{0\}$ is a non-zero central projection, then $p$ is the identity $1$ of $A$ if and only if $\{a\}^\de \cup \{p\}^\de \neq \RS_{A^+}$  for every $a\in \RS_{A^+}$.
\end{lem}

\smnoind

Denote by $\BP_{A^{**}}^\m$ the set of all minimal projections in $A^{**}$, and
\begin{equation}\label{eqt:def-id-atom-bidual}
\ida{A}:= \bigvee \BP_{A^{**}}^\m.
\end{equation}
Then $\ida{A}$ is a central projection in $A^{**}$.
The von Neumann subalgebra
$$
A^{**}_\ba:= \ida{A}\cdd A^{**}
$$
is called the \emph{atomic part} of $A^{**}$, and its predual can be identified with  $\ida{A}\cdd A^*$.

\smnoind

The  \emph{support projection} $\ks_\omega$ of a state $\omega$ on  $A$ is defined by
\begin{equation}\label{eqt:def-s-omega}
	\ks_\omega:= \min \{p\in \BP_{A^{**}}: \omega(p) =1 \}.
\end{equation}
For two states $\omega$ and $\omega'$ on $A$, one has (see, e.g., \cite[Theorem III.4.2(ii)]{Take})
\begin{equation}\label{eqt:pure-st-orth-supp}
	\| \omega - \omega' \|=2\quad  \text{if and only if} \quad \ks_{\omega} \ks_{\omega'} = 0.	
\end{equation}
The assignment
$\omega\mapsto \ks_\omega$ is a bijection from $\ps(A)$ onto $\BP_{A^{**}}^\m$; i.e.,
\begin{equation}\label{eqt:min-pure-st}
	\BP_{A^{**}}^\m = \{\ks_\omega: \omega\in \ps(A) \}.
\end{equation}
In fact, one has
\begin{equation}\label{eqt:omega-supp}
	\omega(a)\ks_\omega = \ks_\omega a \ks_\omega \qquad (\omega\in \ps(A)), a\in A^{**}).
\end{equation}

\smnoind

The following fact will be used implicitly.

\begin{rem}\label{rem:inf-closed}
If $E\subseteq \BP_{A^{**}_\ba}$, then it is easy to see that $\bigwedge E, \bigvee E \in \BP_{A^{**}_\ba}$ and hence, $\bigwedge E$ and $\bigvee E$ coincide with the infimum and the supremum of $E$ in $\BP_{A^{**}_\ba}$.
\end{rem}

Let
$\Phi_A:A^{**}\to A^{**}_\ba$  be the normal $^*$-homomorphism  given by
\begin{equation}\label{eqt:def-Phi-A}
\Phi_A(a) := a\cdd \ida{A} \qquad (a\in A^{**}).
\end{equation}
It is obvious that
\begin{equation}\label{Phi-A-min}
\Phi_A(a) \leq a \ \ \ (a\in A^{**})\qquad \text{and} \qquad \Phi_A(p) = p \ \ \ (p\in \BP_{A^{**}_\ba}),
\end{equation}
which implies
\begin{equation}\label{eqt:Phi-A-dom}
p \leq a \quad \text{if and only if} \quad p \leq \Phi_A(a)\qquad (a\in A^{**}; p\in \BP_{A^{**}_\ba}).
\end{equation}

If $\widehat{A}$ is the set of equivalence classes of irreducible $^*$-representations of $A$, then
\begin{equation}\label{eqt:atom-bidual}
	A^{**}_\ba \cong {\bigoplus}_{(\CH_\pi,\pi)\in \widehat{A}}\, \CL(\CH_\pi),
\end{equation}
and $\Phi_A|_A$ can be identified with the injective $^*$-representation ${\bigoplus}_{(\CH_\pi,\pi)\in \widehat{A}}\ \pi$.

\smnoind

A very important fact (which will be used many times in this article) is that the lattice $\BP_{A^{**}_\ba}$ is \emph{atomistic}, in the sense that
 for each $p\in \BP_{A^{**}_\ba}$ we have
\begin{align}\label{eqt:atomic}
	p = \bigvee \BP_{A^{**}}^\m[p], \quad \text{where} \quad \BP_{A^{**}}^\m[p] :=\big\{e\in \BP_{A^{**}}^\m: e\leq p\big\}.
\end{align}

\begin{rem}\label{rem:proj-atomic}
Consider  $p,q\in \BP_{A^{**}_\ba}$.
The following are consequences of \eqref{eqt:atomic}:
\begin{enumerate}[(i)]
	\item $p\leq q$ if and only if $\BP_{A^{**}}^\m[p]\subseteq \BP_{A^{**}}^\m[q]$.
	
	\item $p\bot q$ if and only if for $e_1,e_2\in \BP_{A^{**}}^\m$ with $e_1\leq p$ and $e_2\leq q$, one has $e_1\bot e_2$.
		
	\item If $p\lneqq q$, then there exists $e\in \BP_{A^{**}}^\m$ such that $e\leq q$ and $e \bot p$.
\end{enumerate}
\end{rem}

As in Definition \ref{defn:open}, we denote by  $\CO(A)$ and $\CC(A)$, respectively, the set of open projections and the set of closed projections of $A$.
As in Definition \ref{defn:q-open}, we denote by  $\CO^\ba(A)$ and $\CC^\ba(A)$, respectively, the set of q-open projections and the set of q-closed projections of $A$.
We recall from \eqref{eqt:Phi-A-closed} that
$$\CO^\ba(A) = \Phi_A(\CO(A))\quad \text{and} \quad \CC^\ba(A) = \Phi_A(\CC(A)).$$
More information on open projections and closed projections can be found in Appendix \ref{sec:open-proj}.

\smnoind


The following is a form  of Dye's theorem for closed projections of $C^*$-algebras.
Although it is a reformulation of \cite[Corollary 3.8]{DN}, we give a self-contained proof here for completeness (see also the proof of \cite[Proposition 2.11]{Ng-spec}).

\begin{prop}[{\cite[Corollary 3.8]{DN}}]\label{prop:DN-3.8}
Let $A$ and $B$ be unital $C^*$-algebras such that $A$ has no 2-dimensional irreducible $^*$-representation.
If $\Gamma: \CC(A) \to \CC(B)$ is a bi-orthogonality preserving bijection, then there is a Jordan $^*$-isomorphism $\Psi:A\to B$ such that $\Gamma = \Psi^{**}|_{\CC(A)}$.
\end{prop}
\begin{proof}
Parts (a) and (b) of Lemma \ref{lem:Phi-A-closed} produce a bijection
$$\Gamma_\ba: \CC^\ba(A) \to \CC^\ba(B)$$
which is a bi-orthogonality preserving map satisfying
\begin{equation}\label{eqt:def-Gamma-a}
\Gamma_\ba\circ \Phi_A|_{\CC(A)} = \Phi_B\circ \Gamma.
\end{equation}
As $\BP_{A^{**}}^\m \subseteq \CC^\ba(A)$ (see Lemma \ref{lem:Phi-A-closed}(c)), we may define $\ti \Gamma_\ba:\BP_{A^{**}_\ba}\to \BP_{B^{**}_\ba}$ by
$$\ti \Gamma_\ba(p):= \bigvee \Gamma_\ba\big(\BP_{A^{**}}^\m[p]\big)\qquad (p\in \BP_{A^{**}_\ba}).$$

We first claim that $\Gamma_\ba$ is an order isomorphism.
In fact,
since $1_{A^{**}_\ba}$ is the only element in $\CC^\ba(A)$ admitting just one element in $\CC^\ba(A)$ orthogonal to it (namely, $0$),
$$\Gamma_\ba(1_{A^{**}_\ba}) = 1_{B^{**}_\ba}.$$
Let $p,q\in \CC^\ba(A)\setminus \{1_{A^{**}_\ba}\}$ with $p\leq q$.
If $f\in \BP_{B^{**}}^\m[1_{B^{**}_\ba}- \Gamma_\ba(q)]$, then we learn from $f\bot \Gamma_\ba(q)$ that $\Gamma_\ba^{-1}(f)\bot q$, which gives $\Gamma_\ba^{-1}(f)\bot p$, and hence $f\in \BP_{B^{**}}^\m[1_{B^{**}_\ba}- \Gamma_\ba(p)]$.
Thus, Remark \ref{rem:proj-atomic}(i) implies that $1_{B^{**}_\ba}-\Gamma_\ba(q) \leq 1_{B^{**}_\ba}-\Gamma_\ba(p)$.
This means that $\Gamma_\ba$ is order preserving.
In the same way, $\Gamma_\ba^{-1}$ is also order preserving.

Since $\BP_{A^{**}}^\m$ is the set of minimal elements in $\CC^\ba(A)\setminus \{0\}$,
the above claim produces
\begin{equation}\label{eqt:Gamma-0-min-proj}
\Gamma_\ba(\BP_{A^{**}}^\m) = \BP_{B^{**}}^\m.
\end{equation}
as well as
$$\Gamma_\ba\big(\BP_{A^{**}}^\m[p]\big) = \BP_{B^{**}}^\m\big([\Gamma_\ba(p)]\big) \qquad (p\in \CC^\ba(A)).$$
Thus, \eqref{eqt:atomic} implies that
\begin{equation}\label{eqt:ti-Gamma-0}
\tilde \Gamma_\ba (p)  =\Gamma_\ba (p) \qquad (p\in \CC^\ba(A)).
\end{equation}

On the other hand, it follows from the definition that $\ti \Gamma_\ba$ is order preserving.
Moreover, the definition of $\ti \Gamma_\ba$, \eqref{eqt:atomic}, Remark \ref{rem:proj-atomic}(ii), \eqref{eqt:Gamma-0-min-proj} and the fact that $\Gamma_\ba$ is bi-orthogonality preserving, tell us that for each $p\in \BP_{A^{**}_\ba}$,
\begin{align}
	\ti \Gamma_\ba(1_{A^{**}_\ba}-p)
	&  = \bigvee \big\{\Gamma_\ba(e):  e\in \BP_{A^{**}}^\m \text{ and } e\bot p\big\}\nonumber\\
	&  = \bigvee \big\{\Gamma_\ba(e):  e\in \BP_{A^{**}}^\m \text{ such that } e\bot e' \text{ for each }e'\in \BP_{A^{**}}^\m[p]\big\}\nonumber\\
	&  =  \bigvee \big\{f\in \BP_{B^{**}}^\m: f\bot f' \text{ for each }f'\in \Gamma_\ba\big(\BP_{A^{**}}^\m[p]\big)\big\}\nonumber\\
	& = \bigvee \big\{f\in \BP_{B^{**}}^\m: f\bot \tilde\Gamma_\ba(p)\big\}
	= 1_{B^{**}_\ba}- \tilde\Gamma_\ba(p). \label{eqt:pres-ortho-comp}
\end{align}

Similarly, we obtain a map $\ti \Gamma_\ba^\flat:\BP_{B^{**}_\ba}\to \BP_{A^{**}_\ba}$ such that  $\ti \Gamma_\ba^\flat(r):= \bigvee \Gamma_\ba^{-1}\big(\BP_{B^{**}}^\m[r]\big)$ for each $r\in \BP_{B^{**}_\ba}$.
As in the above, $\ti \Gamma_\ba^\flat$ is order preserving and satisfies
\begin{equation*}
\ti \Gamma_\ba^\flat(1_{B^{**}_\ba}-r) =1_{A^{**}_\ba} - \ti \Gamma_\ba^\flat(r).
\end{equation*}
Let $p\in \BP_{A^{**}_\ba}$.
One clearly has
$$\BP_{A^{**}}^\m[p] \subseteq \Gamma_\ba^{-1}\big(\BP_{B^{**}}^\m\big(\ti \Gamma_\ba(p)\big)\big),$$
which gives $p \leq \ti\Gamma_\ba^\flat\big(\ti\Gamma_\ba(p)\big)$ (see \eqref{eqt:atomic}).
In the same way, we also have
$$
 1_{A^{**}_\ba}- p \leq \ti\Gamma_\ba^\flat\big(\ti\Gamma_\ba(1_{A^{**}_\ba}-p)\big).
 $$
Thus,
$$p \geq 1_{A^{**}_\ba} - \ti\Gamma_\ba^\flat\big(\ti\Gamma_\ba(1_{A^{**}_\ba}-p)\big) = \ti\Gamma_\ba^\flat\big(1_{B^{**}_\ba} - \ti\Gamma_\ba(1_{A^{**}_\ba}-p)\big) = \ti\Gamma_\ba^\flat\big(\ti\Gamma_\ba(p)\big).$$
Therefore, $p = \ti\Gamma_\ba^\flat\big(\ti\Gamma_\ba(p)\big)$.
By symmetry, $r = \ti\Gamma_\ba\big(\ti\Gamma_\ba^\flat(r)\big)$ ($r\in \BP_{B^{**}_\ba}$).
This means that $\ti \Gamma_\ba$ is an order isomorphism.
From this and \eqref{eqt:pres-ortho-comp}, we know that $\ti \Gamma_\ba$ is a bi-orthogonality preserving bijection.

Now, Dye's theorem (see \cite{dye1955geometry}) tells us that $\ti\Gamma_\ba$ extends to a Jordan $^*$-isomorphism $\Theta:A^{**}_\ba \to B^{**}_\ba$.
In particular, \eqref{eqt:ti-Gamma-0} ensures that
$$\Theta|_{\CC^\ba(A)} = \Gamma_\ba, \quad \text{and} \quad \Theta(\CC^\ba(A)) = \CC^\ba(B).$$
It then follows from \cite[Theorem I.1]{Ake71} that $\Theta(\Phi_A(A)) = \Phi_B(B)$, and we obtain from the injectivity of $\Phi_B$ a Jordan $^*$-isomorphism $\Psi:A\to B$ such that
$$\Theta\circ \Phi_A|_A = \Phi_B\circ \Psi.$$
This equality implies that $\Theta \circ \Phi_A = \Phi_B\circ \Psi^{**}$, and hence (see \eqref{eqt:def-Gamma-a})
$$\Phi_B\circ \Gamma = \Gamma_\ba\circ \Phi_A|_{\CC(A)} = \Theta \circ \Phi_A|_{\CC(A)} = \Phi_B\circ \Psi^{**}|_{\CC(A)}.$$
It then follows from Lemma \ref{lem:Phi-A-closed}(a) that $\Gamma = \Psi^{**}|_{\CC(A)}$.
\end{proof}

Note that the argument  above also produces a proof for a version of
Dye's theorem for q-closed projections of $C^*$-algebras (see \cite[Theorem 3.6]{DN}).

\begin{defn}\label{def:peak-zero-supp}
Let $A$ be a unital $C^*$-algebra.
Let $u\in \RS_{A^+}$ and $a\in \RB_{A^+}$.

\smnoind
(a) The \emph{peak projection} $\bp_A(u)$ and the \emph{atomic peak projection} $\bp_A^\ba(u)$ of $u$
are
\begin{equation*}\label{eqt:def-peak}
\bp_A(u) := w^*\text{-}{\lim}_n\, u^n \in \CC(A) \quad \text{and} \quad \bp_A^\ba(u):= \Phi_A(\bp_A(u))\in \CC^\ba(A).
\end{equation*}

\smnoind
(b) The \emph{null projection} $\bz_A(a)$ and the \emph{atomic null projection} $\bz_A^\ba(a)$ of $a$ are
\begin{equation*}\label{eqt:def-null}
\bz_A(a) := w^*\text{-}{\lim}_n\, (1-a)^n \in \CC(A)	\quad \text{and} \quad \bz_A^\ba(a):= \Phi_A(\bz_A(a))\in \CC^\ba(A).
\end{equation*}

\smnoind
(c) The \emph{support projection} $\bs_A(a)$ and the \emph{atomic support projection} $\bs_A^\ba(a)$ of $a$ are
\begin{equation*}
\bs_A(a) = w^*\text{-}{\lim}_n\, a^{1/n}\in \CO(A) \quad \text{and} \quad \bs_A^\ba(a):= \Phi_A(\bs_A(a))\in \CO^\ba(A).
\end{equation*}
\end{defn}

Notice that $\bp_A(u), \bz_A(a) \in \CC(A)$, and  $\bs_A(a) \in \CO(A)$ because of
Remark \ref{rem:meet-semilat}(a)
Moreover, $\bz_A(v) = 0$ for every $v\in \RS_{A^+}^{-1}$ and we have (see \eqref{eqt:atomic})
\begin{equation}\label{eqt:peak-non-zero}
\bp_A^\ba(u)\neq 0,\quad \text{equivalently,} \quad \BP_{A^{**}}^\m[\bp_A^\ba(u)]\neq \emptyset, \quad  \text{for any} \ u\in \RS_{A^+},
\end{equation}
because there is $e\in \BP_{A^{**}}^\m$ with $eue  =e$ (see \eqref{eqt:omega-supp}) and hence $eu = e$.

\smnoind

Suppose that $B\subseteq A$ is a unital $C^*$-subalgebra of $A$.
One may identify $B^{**}$ as a von Neumann subalgebra of $A^{**}$ in the canonical way.
Under this identification,
\begin{equation}\label{eqt:equal-peak-supp}
	\bp_B(u) = \bp_A(u), \quad \bs_B(a) = \bs_A(a)  \quad \text{and}\quad \bz_B(a) = \bz_A(a),
	\end{equation}
where $u\in \RS_{B^+}$ and $a\in \RB_{B^+}$.
However, it is possible to have
$\bp^\ba_B(u)\neq \bp^\ba_A(u)$, $\bz^\ba_B(u) \neq \bz^\ba_A(u)$ or $\bs^\ba_B(u) \neq \bz^\ba_A(u).$
For a counter example, consider the inclusion $C[0,1]\subseteq \CL(L^2[0,1])$.



\begin{lem}\label{lem:peak-zero}
	Let $A$ be a unital $C^*$-algebra.
	Consider $u\in \RS_{A^+}$, $a\in \RB_{A^+}$, $s\in \BP_{A^{**}}$ and $p\in \BP_{A^{**}_\ba}$.
	
\smnoind
(a) $\bp_A(u) = \bz_A(1-u)$.

\smnoind
(b) $\bs_A(a) = 1 -\bz_A(a)$.

\smnoind
(c) $\bs_A(a) = \min\{r\in \BP_{A^{**}}: ar = a\}$.
Hence, $\bs_A(a)\leq s$ if and only if $a\leq s$.

\smnoind
(d) $\bz_A(a) = \max\{r\in \BP_{A^{**}}:  ar = 0\}$ as well as $\bz_A^\ba(a) = \bigvee \{e\in \BP_{A^{**}}^\m:  ae = 0\}$.
Furthermore, $s\leq \bz_A(a)$  if and only if $s\leq 1-a$, and $p\leq \bz_A^\ba(a)$ if and only if $p\leq \Phi_A(1-a)$.

	\smnoind
(e) $\bp_A(u) = \max\{r\in \BP_{A^{**}}:  ur = r\}$ and $\bp_A^\ba(u) = \bigvee \{e\in \BP_{A^{**}}^\m:  ue = e\}$.
Moreover, $s\leq \bp_A(u)$ if and only if $s\leq u$, and $p\leq \bp_A^\ba(u)$ if and only if $p\leq \Phi_A(u)$.

	\smnoind
	(f) $\bp_A(u) \bz_A(u) = 0$, and hence, $\bp_A^\ba(u) \bz_A^\ba (u) = 0$.
	
	\smnoind
	(g) If $\bp_A(u) + \bz_A(u) =1$, then $u\in \BP_A$.
\end{lem}
\begin{proof}
(a) This part follows from the definition.
	
\smnoind
(b) Let $B\subseteq A$ be the $C^*$-subalgebra generated by $a$ and $1$.
Then $B \cong C(\sigma(a))$ and $B^{**}_\ba\cong \ell^\infty(\sigma(a))$.
Under this identification, we have
$$\bs^\ba_B(a) = w^*\text{-}{\lim}_n\, \Phi_B(a)^{1/n} = \{x\in \sigma(a):a(x) > 0\}$$
and $\bz^\ba_B(a) = \{x\in \sigma(a):a(x) = 0\}$.
This means that $\bs^\ba_B(a) = 1_{B^{**}_\ba} -\bz^\ba_B(a)$.
Since $\bs_B(a), 1-\bz_B(a)\in \CO(B)$, we learn from Lemma \ref{lem:Phi-A-closed}(a) that
$\bs_B(a) = 1 -\bz_B(a).$
Thus, \eqref{eqt:equal-peak-supp} gives
$$\bs_A(a) = 1 -\bz_A(a).$$

	\smnoind
	(c) Let $r\in \BP_{A^{**}}$ such that $a r = a$.
	As $ar=ra$ as well, we know that $a^{1/n}r =a^{1/n}$.
	This gives $\bs_A(a)r = \bs_A(a)$, i.e., $\bs_A(a)\leq r$.
	On the other hand, it follows from $\|aa^{1/n} -a\|\to 0$ that $a \bs_A(a) = a$.
Hence, $\bs_A(a) =  \min\{r\in \BP_{A^{**}}: ar = a\}$, and the first claim is established.
The second claim follows from the first claim.

\smnoind
(d) The first equality follows from parts (b) and (c) above.
Let $e\in \BP_{A^{**}}^\m$.
If $e \leq \bz_A^\ba(a)$, then $e\leq \bz_A(a)$ (as $\bz_A^\ba(a) = \bz_A(a)\cdd \ida{A}$) and the first equality gives $ae =0$.
Conversely, if $ae = 0$, then the first equality implies that $e\leq \bz_A(a)$, and so, $e = \Phi_A(e)\leq \bz_A^\ba(a)$ (see   \eqref{Phi-A-min}).
This means that
$$\BP_{A^{**}}^\m(\bz_A^\ba(a)) = \{e\in \BP_{A^{**}}^\m: ae=0 \}$$
and the second equality follows from \eqref{eqt:atomic}.

Furthermore, we learn from parts (b) and (c) above that $r\leq \bz_A(a)$  if and only if $r\leq 1-a$.
This, together with  \eqref{Phi-A-min} and \eqref{eqt:Phi-A-dom}, tells us that $p\leq \bz_A^\ba(a)$ if and only if $p\leq \Phi_A(1-a)$.

	\smnoind
	(e) This part follows from parts (a) and (d) above.
	
	\smnoind
	(f) By part (d) above, one has $u^m\bz_A(u) = 0$ ($m\in \BN$), which gives the conclusion.
	
	\smnoind
	(g) 
Let $D\subseteq A$ be the $C^*$-subalgebra generated by $u$ and $1$.
Then $D \cong C(\sigma(u))$.
	As in \eqref{eqt:equal-peak-supp}, we have $\bp_A(u) = \bp_D(u)$ and $\bz_A(u) = \bz_D(u)$.
	It then follows from the hypothesis that $\bp_D(u) + \bz_D(u) = 1$, and hence $\bp_D^\ba(u) + \bz_D^\ba(u) = 1_{D^{**}_\ba}$.
	Under the identification of $D \cong C(\sigma(u))$ and $D^{**}_\ba \cong \ell^\infty(\sigma(u))$, we have
	$$\bz^\ba_D(u) = w^*\text{-}{\lim}_n \big(1-\Phi_D(u)\big)^n = \{x\in \sigma(u):u(x) = 0\},$$
	and $\bp^\ba_D(u)  = \{x\in \sigma(u):u(x) =1 \}$.
	We then see from $\bp_D^\ba(u) + \bz_D^\ba(u) = 1_{D^{**}_\ba}$ that  $\Phi_D(u)\in \ell^\infty(\sigma(u))$ is a projection, and hence so is $u$ (because $\Phi_D|_D$ is injective).
\end{proof}

\begin{prop}\label{prop:dist-one}
	Let $A$ be a unital $C^*$-algebra.
	Let  $\omega\in \ps(A)$, $u,w\in \RS_{A^+}$ and $a\in \RB_{A^+}$.
	Denote by $\ks_\omega$ the support projection of $\omega$ (see \eqref{eqt:def-s-omega}).
	
	\smnoind
	(a) $\omega(a) = 0$ if and only if $\ks_\omega\leq \bz_A(a)$; equivalently, $\ks_\omega\leq \bz_A^\ba(a)$. 	
		
	\smnoind
	(b) $\omega(u) = 1$ if and only if $\ks_\omega\leq \bp_A(u)$; equivalently, $\ks_\omega\leq \bp_A^\ba(u)$.

	\smnoind
	(c) $\|u-w\| =1$ if and only if either $\bp_A(u) \wedge \bz_A(w) \neq 0$ or $\bz_A(u) \wedge \bp_A(w) \neq 0$; equivalently, either
 $\bp_A^\ba(u) \wedge \bz_A^\ba(w) \neq 0$ or $\bz_A^\ba(u) \wedge \bp_A^\ba(w) \neq 0$.
	
	\smnoind
	(d) $w\in \RS_{A^+}^{-1}$ if and only if $\bz_A(w) = 0$; equivalently, $\bz_A^\ba(w) = 0$.
\end{prop}
\begin{proof}
By \eqref{Phi-A-min}, parts (a) and (b) of Lemma \ref{lem:Phi-A-closed} as well as Remark \ref{rem:inf-closed=closed}(b), we know that for every $e\in \BP_{A^{**}}^\m$ and $r,s\in \CC(A)$, one has $r\wedge s\in \CC(A)$, and
$$e\leq r\quad  \text{if and only if} \quad e\leq \Phi_A(r)$$
as well as $r\wedge s =0$ if and only if $\Phi_A(r)\wedge \Phi_A(s) = 0$.
Thus, it suffices to verify the corresponding statements concerning $\bp_A^\ba$ and $\bz_A^\ba$. 	

\smnoind
(a) Suppose that $\omega(a) = 0$.
As $\ks_\omega a \ks_\omega = \omega(a)\ks_\omega = 0$ (see \eqref{eqt:omega-supp}), we have $a \ks_\omega = 0$, and hence
$\ks_\omega  \leq \bz_A^\ba(a)$ (see Lemma \ref{lem:peak-zero}(d)).
Conversely, if $\ks_\omega\leq \bz_A^\ba(a)$,  it follows from $\bz_A^\ba(a) \leq 1-a$ that
$$ 1 = \omega(\ks_\omega) \leq \omega(1-a) \leq 1,$$
and hence $\omega(a) = 0$.

\smnoind
(b)	This follows from part (a) above and Lemma \ref{lem:peak-zero}(a).
	
	\smnoind
	(c) 	Suppose that $\|u-w\| =1$.
	By Lemma \ref{lem:peak-and-zero}(a), there exists $\rho\in \ps(A)$ such that $\{\rho(u), \rho(w)\} = \{0,1\}$.
	Therefore, by parts (a) and (b) above, we either have $\ks_\rho \leq \bp_A^\ba(u) \wedge \bz_A^\ba(w)$ or $\ks_\rho \leq \bz_A^\ba(u) \wedge \bp_A^\ba(w)$.
	
	Conversely, suppose that $\bp^\ba_A(u) \wedge \bz^\ba_A(w)\neq 0$.
	One can find $e\in \BP_{A^{**}}^\m$ such that  $e\leq \bp^\ba_A(u) \wedge \bz^\ba_A(w)$.
	Relation \eqref{eqt:min-pure-st} then produces  $\omega_e\in \ps(A)$ satisfying $\ks_{\omega_e} = e$.
	Now, parts (a) and (b) above tell us that $\omega_e(u) = 1$ and $\omega_e(w) = 0$, which gives $\|u-w\| =1$.
	Similarly, if $\bz_A^\ba(u) \wedge \bp_A^\ba(w)\neq 0$, then $\|u-w\| =1$.
	
	\smnoind
	(d) This follows directly from Lemma \ref{lem:peak-and-zero}(b) as well as part (a) above.
\end{proof}


\smnoind

\begin{rem}\label{rem:q-closed}
Let $a\in \RB_{A^+}$.
Then $\overline{aAa}$ is the smallest hereditary $C^*$-subalgebra of $A$ containing $a$, and  $(\overline{aAa})^{**}$ can be identified with the smallest weak$^*$-closed hereditary $C^*$-subalgebra of $A^{**}$ containing $a$.
This and Lemma \ref{lem:peak-zero}(c) give
	\begin{align*}
		\overline{aAa} = \bs_A(a)\cdd A^{**}\cdd \bs_A(a)\cap A.
	\end{align*}
\end{rem}

The following is taken from \cite{peralta2018characterizing}.
We give a slightly different proof here using Proposition \ref{prop:dist-one}(c).

\begin{cor}[{\cite[Proposition 2.2]{peralta2018characterizing}}] \label{cor:Peral-2.2}
	Suppose that $A$ is a unital $C^*$-algebra.
	If $u\in \RS_{A^+}$ satisfying  $\{u\} = (\{u\}^\de)^\de$, then $u$ is a projection.
\end{cor}
\begin{proof}
	Let  $B\subseteq A$ be the $C^*$-algebra generated by $u$ and $1$.
	One has  $B\cong C(\sigma(u))$.
	Assume, on the contrary, that $\sigma(u)\not\subseteq \{0,1\}$.
	Define $w\in B$ by $w(t):=t^2$ for each $t\in \sigma(u)$.
	Then  $w(0) = 0$ if $0\in \sigma(u)$ and  $w(1) = 1$.
	By the assumption, $w\neq u$.
However, under the identification of $B^{**}_\ba\cong \ell^\infty(\sigma(u))$, one has
$$\bp_B^\ba(w) = w^*\text{-}{\lim}_n\, \Phi_B(w)^n = \{x\in \sigma(u):w(x) =1 \} = \bp_B^\ba(u)$$
and $\bz_B^\ba(w) = \{x\in \sigma(u):w(x) = 0 \} = \bz_B^\ba(u)$.
Thus, \eqref{eqt:equal-peak-supp} and Lemma \ref{lem:Phi-A-closed}(a) give
$$\bp_A(w) = \bp_B(w) = \bp_B(u) = \bp_A(u)$$
and $\bz_A(w) = \bz_A(u)$.
	Proposition \ref{prop:dist-one}(c) then implies that $\{w\}^\de = \{u\}^\de$.
	Hence, $w\in (\{u\}^\de)^\de$, but the hypothesis will then produce the contradiction that $w = u$.
\end{proof}

The converse of Corollary \ref{cor:Peral-2.2} will be given
in Corollary \ref{cor:met-char-proj} below.

\section{Order and algebraic properties that determined by diametral relation}\label{S:projections}

In this section, we will present some order properties and algebraic properties that are determined by the diametral relation $\D_{\RS_{A^+}}$ on $\RS_{A^+}$ (see \eqref{eqt:def-diam-rel}).
Obviously, information of $\D_{\RS_{A^+}}$ are captured by the assignment (see \eqref{eqt:def-orth-comp})
$$u\mapsto \{u\}^\de = \{v\in \RS_{A^+}: \|v-u\| =1 \}$$
and we will  concern with $\{u\}^\de$ for $u\in \RS_{A^+}$.

\smnoind

Let us begin this section with the following lemma, which is a generalization of \cite[Lemma 2.3(c)]{LN} and \cite[Lemma 4.7]{LNW}.

\begin{lem}\label{lem:p-leq-a}
	Let $A$ be a unital $C^*$-algebra.
	For $q\in \BP_A\setminus \{0,1\}$ and $a\in \RS_{A^+}^{-1}$, one has $q\leq a$ if and only if $1-b \in \{a\}^\de$ for every
	$b\in \RS_{A^+}$ with $b\leq q$.
\end{lem}
\begin{proof}
	The necessity follows from the fact that when $q\leq a$, one has $1-a\leq 1-q$, which implies $\|(1-b) -a \| = \|(1-a) - b\| = \max\{\|1-a\|, \|b\| \} = 1$.
	
	For the sufficiency, one needs to verify that $q(1-a) =0$, when the said condition in the statement holds.
	Assume, on the contrary, that $\kappa:= \|q(1-a)q\| > 0$.
	Set
	$$b:= \frac{q(1-a)q}{\kappa}.$$
	Then $b\in \RS_{A^+}$ and $b\leq q$.
	Let
	$$s:= \bs_A(b)\in \CO(A).$$
	By the hypothesis,  we have $\|(1-b)-a\| = 1$.
	Thus, by parts (a) and (b) of Lemma \ref{lem:peak-and-zero}, there exists $\omega\in \ps(A)$ satisfying
	\begin{equation}\label{eqt:omega-a-1-b}
		\omega(a) = 1 \quad \text{and} \quad \omega(1-b) = 0.
	\end{equation}
	The second equality above, together with Proposition \ref{prop:dist-one}(a) as well as parts (a), (f) and (b) of Lemma \ref{lem:peak-zero},  implies that (see \eqref{eqt:def-s-omega})
	$$\ks_\omega \leq \bz_A(1-b) = \bp_A(b) \leq 1 - \bz_A(b) = s.$$
	Thus, we may use \eqref{eqt:omega-supp} to get
	$$\omega(s) =1.$$
	On the other hand, as $b\leq q$, Lemma \ref{lem:peak-zero}(c) implies $\bs_A(b)\leq q$; i.e., $sq =s$. Thus,
	$$\kappa b = s(\kappa b)s = s(q(1-a)q)s = s(1-a)s$$
	(Lemma \ref{lem:peak-zero}(c) gives the first equality above).
	From this, \eqref{eqt:omega-a-1-b} and Lemma \ref{lem:peak-and-zero}(c), 
	$$\kappa = \kappa \omega (b) = \omega(s(1-a)s) = \omega(1-a) = 1 - \omega(a) = 0,$$
	which is a contradiction.
\end{proof}

Our next lemma is an application of \cite[Theorem 1.1]{Ake} and could be a known fact, but we give its argument here for completeness.

\begin{lem}\label{lem:one-peak}
Let $B$ be a separable unital $C^*$-algebra.

\smnoind
(a) If $\omega_0,\omega_1\in \ps(B)$ satisfying $\ks_{\omega_0}\ks_{\omega_1} = 0$, then one can find $a\in \RS_{B^+}$ such that
\begin{equation}\label{eqt:one-peak}
\omega_0(a) =0, \quad \omega_1(a) =1 \quad \text{and}\quad 0 < \omega'(a) < 1 \ \  (\omega'\in \ps(B)\setminus \{\omega_0,\omega_1\}).
\end{equation}

\smnoind
(b) For any $\omega_1\in \ps(B)$, there exists $b\in \RS_{B^+}^{-1}$ such that
\begin{equation}\label{eqt:one-peak-at-omega-e}
	\omega_1(b) =1 \quad \text{and}\quad 0<\omega'(b) < 1 \ \  (\omega'\in \ps(B)\setminus \{\omega_1\}).
\end{equation}
\end{lem}
\begin{proof}
(a)  Since $\|\omega_0 - \omega_1\| =2$ (see \eqref{eqt:pure-st-orth-supp}), one obtains from \cite[Theorem 1.1]{Ake} a maximal abelian $C^*$-subalgebra $D\subseteq B$ such that $\nu_k:=\omega_k|_D\in \ps(D)$, and
$$
\omega_k \text{ is the only element in $\ps(B)$ extending $\nu_k$ for $k=0,1$.}
$$
Since $D$ is a separable unital abelian $C^*$-algebra, there is a compact metrizable space $X$ such that $D \cong C(X)$.
In this case, $\nu_k$ is the point mass at some $x_k\in X$, for $k=0,1$.
As $X$ is compact metrizable, there exists $a\in \RS_{D^+}\subseteq \RS_{B^+}$ such that
\begin{equation}\label{eqt:one-peak-comm}
\big\{y\in X: a(y) = 0 \big\} = \{x_0\} \quad \text{and}\quad \big\{y\in X: a(y) =1 \big\} = \{x_1\}.
\end{equation}
Let $\omega'\in \ps(B)$ with $\omega'(a) = 0$.
Since $\omega'|_D$ is a state, there exists a probability measure $\mu'$ on $X$ satisfying
$\omega'(c) = \int c\, d\mu'$ ($c\in D$).
In particular,
$$1 = \omega'(1-a) = \int (1-a)\, d\mu'.$$
From this, the support of $\mu'$ is contained in $\{y\in X: (1-a)(y) =1 \} = \{x_0\}$ (see \eqref{eqt:one-peak-comm}).
This means that $\mu'$ is the point mass at $x_0$; i.e. $\mu' = \nu_0$.
Since $\omega_0$ is the unique extension of $\nu_0$ to $B$,
we see that
$\omega' = \omega_0$.
In the same way, if $\omega'\in \ps(B)$ satisfying $\omega'(a) =1$, then $\omega' = \omega_1$.

\smnoind
(b) The statement is clear when $B\cong \BC$.
Assume that $B\neq \BC$.
Pick any $\omega_0\in \ps(B)$ with $\ks_{\omega_0}\ks_{\omega_1} = 0$.
Let $a\in \RS_{B^+}$ be the element as given by part (a) above.  Then
$$
b:= (1+a)/2\in \RS_{B^+}^{-1}
$$
will satisfy the requirement.
\end{proof}

The above gives the following non-commutative analogue of Urysohn's Lemma.
Note that this analogue  is different from \cite[Theorem I.1]{Ake70}.
In fact, unlike \cite[Theorem I.1]{Ake70}, the following result  does not hold without the separability assumption, even in the commutative case.

\begin{cor}\label{cor:Ury}
Let $B$ be a separable unital $C^*$-algebra.

\smnoind
(a) For $e_0,e_1\in \BP_{B^{**}}^\m$ with $e_0e_1=0$, there exists $a\in \RS_{B^+}$ such that $\bp_B(a) =  \bp_B^\ba(a) = e_1$ and $\bz_B(a) = \bz_B^\ba(a) = e_0$.
	
\smnoind
(b) For $e_1\in \BP_{B^{**}}^\m$, there exists $b\in \RS_{B^+}^{-1}$ such that $\bp_B(b) = \bp_B^\ba(b) = e_1$.
\end{cor}
\begin{proof}
(a) By \eqref{eqt:min-pure-st}, there are $\omega_0,\omega_1\in \ps(B)$ satisfying $\ks_{\omega_i} = e_i$ ($i=1,2$).
Let $a\in \RS_{B^+}$ be the element obtained in Lemma \ref{lem:one-peak}(a).
Since $\omega_0(a) =0$, $\omega_1(a) =1$ and $0 < \omega'(a) < 1$ for each $\omega'\in \ps(B)\setminus \{\omega_0,\omega_1\}$, we learn from parts (a) and (b) of Proposition \ref{prop:dist-one} that $e_0\leq \bz_B^\ba(a)$, $e_1\leq \bp_B^\ba(a)$, as well as
$$e\nleq \bz_B^\ba(a) \quad \text{and} \quad e\nleq \bp_B^\ba(a) \qquad (e\in \BP_{A^{**}}^\m \setminus \{e_0,e_1\}).$$
Hence,  \eqref{eqt:atomic} implies that
$$\bz_B^\ba(a)= e_0 \in \BP_{A^{**}}^\m \quad \text{and} \quad \bp_B^\ba(a) = e_1\in \BP_{A^{**}}^\m.$$
It then follows from \eqref{Phi-A-min}, parts (a) and (b) of Lemma \ref{lem:Phi-A-closed} as well as \eqref{eqt:min>closed} that
$$e_0 = \bz_B(a)\quad \text{and} \quad e_1 = \bp_B(a).$$

\smnoind
(b) The argument is similar to that of part (a) above and will be omitted.
\end{proof}

For a subset $E\subseteq \BP_{A^{**}}^\m$, we set
\begin{equation}\label{eqt:def-bot}
E^\bot := \{f\in \BP_{A^{**}}^\m: ef =0, \text{ for every }e\in E\}.
\end{equation}
It is not hard to see that if $e\in \BP_{A^{**}}^\m$, then $\{e\}^\bot = \BP_{A^{**}}^\m[1-e]$ and
\begin{equation}\label{eqt:e-bot-bot}
(\{e\}^\bot)^\bot = \{e\}.
\end{equation}

\begin{prop}\label{prop:order-1}
Let $A$ be a unital $C^*$-algebra, and $u,v\in \RS_{A^+}$.

\smnoind
(a) $\{u\}^\de \cap \RS_{A^+}^{-1}\subseteq \{v\}^\de \cap \RS_{A^+}^{-1}$ if and only if $\bz_A(u) \leq \bz_A(v)$.

\smnoind
(b) $\{u\}^\de \subseteq \{v\}^\de$ if and only if $\bp_A(u) \leq \bp_A(v)$ and $\bz_A(u) \leq \bz_A(v)$.
\end{prop}
\begin{proof}
Let $B$ be the unital $C^*$-subalgebra of $A$ generated by $u$ and $v$.
Then $B$ is separable.
As in \eqref{eqt:equal-peak-supp}, one has
\begin{equation}\label{eqt:cp-A-B}
\bp_B(u) = \bp_A(u), \quad \bz_B(u) = \bz_A(u), \quad \bp_B(v) = \bp_A(v) \quad \text{and}\quad \bz_B(v) = \bz_A(v).
\end{equation}

\smnoind
(a) By parts (c) and (d) of Proposition \ref{prop:dist-one}, one knows that if $\bz_A(u) \leq \bz_A(v)$  then $\{u\}^\de \cap \RS_{A^+}^{-1} \subseteq \{v\}^\de \cap \RS_{A^+}^{-1}$.

Conversely, assume that $\bz_A(u) \nleq \bz_A(v)$.
By \eqref{eqt:cp-A-B} and Lemma \ref{lem:Phi-A-closed}(b), one has $\bz^\ba_B(u) \nleq \bz^\ba_B(v)$.
We then learn from Remark \ref{rem:proj-atomic}(i) that
$$\big\{e\in \BP_{B^{**}}^\m: e\leq \bz^\ba_B(u) \big\} \nsubseteq \big\{e\in \BP_{B^{**}}^\m: e\leq \bz^\ba_B(v) \big\}.$$
Pick $e_1\in \BP_{B^{**}}^\m$ such that
$$e_1\leq \bz^\ba_B(u)\quad  \text{but} \quad e_1\nleq \bz^\ba_B(v).$$
As $B$ is separable, Corollary \ref{cor:Ury}(b) gives an element $b\in \RS_{B^+}^{-1}$ with $\bp_B^\ba(b) = e_1$.

We claim that $b\notin \{v\}^\de$.
Assume, on the contrary that, $\|b-v\| =1$.
One obtains from parts (c) and (d) of Proposition \ref{prop:dist-one}  as well as \eqref{eqt:atomic} a minimal projection $f\in \BP_{B^{**}}^\m$ such that $f\leq \bp_B^\ba(b)$ and $f\leq \bz_B^\ba(v)$.
As $\bp_B^\ba(b) = e_1\in \BP_{B^{**}}^\m$, we see that $f=e_1$, but this gives the contradiction that $e_1\leq \bz^\ba_B(v)$.

On the other hand, since $e_1\leq \bz^\ba_B(u)$ and $\bp_B^\ba(b) = e_1$, we have $b\in \{u\}^\de$
(see Proposition \ref{prop:dist-one}(c)).
Thus, $\{u\}^\de \cap \RS_{A^+}^{-1} \nsubseteq \{v\}^\de \cap \RS_{A^+}^{-1}$ (note that $b\in \RS_{B^+}^{-1}\subseteq \RS_{A^+}^{-1}$).

\smnoind
(b) Proposition \ref{prop:dist-one}(c) tells us that  if $\bp_A(u) \leq \bp_A(v)$ and $\bz_A(u) \leq \bz_A(v)$, then one has $\{u\}^\de \subseteq \{v\}^\de$.

For the converse, let us assume that $\{u\}^\de \subseteq \{v\}^\de$.

Suppose that $\bz_A(u) \nleq \bz_A(v)$.
Then part (a) above gives  $b\in \RS_{A^+}^{-1}\cap \{u\}^\de\setminus \{v\}^\de$, which contradicts the assumption.

Suppose that $\bp_A(u) \nleq \bp_A(v)$.
In other words, $\bp^\ba_B(u) \nleq \bp^\ba_B(v)$ (see \eqref{eqt:cp-A-B}  and Lemma \ref{lem:Phi-A-closed}(b)).
As in the proof of part (a) above, one can find $e_0\in \BP_{B^{**}}^\m$ with
\begin{equation}\label{eqt:e-0}
e_0\leq \bp^\ba_B(u)\quad  \text{but} \quad e_0\nleq \bp^\ba_B(v).	
\end{equation}

We first claim that (see \eqref{eqt:atomic})
\begin{equation}\label{eqt:e-0-bot}
\{e_0\}^\bot \not\subseteq \BP_{A^{**}}^\m[\bz^\ba_B(v)] =\{f\in \BP_{B^{**}}^\m: f\leq \bz^\ba_B(v)\}.
\end{equation}
Indeed, if
$\{e_0\}^\bot \subseteq \BP_{A^{**}}^\m[\bz^\ba_B(v)]$,
then by \eqref{eqt:e-bot-bot}, Lemma \ref{lem:peak-zero}(f) and \eqref{eqt:peak-non-zero}, one has
$$\{e_0\} = (\{e_0\}^\bot)^\bot \supseteq \big\{f\in \BP_{B^{**}}^\m: f\leq \bz^\ba_B(v)\big\}^\bot \supseteq \{e\in \BP_{B^{**}}^\m: e\leq \bp^\ba_B(v)\}\neq \emptyset,$$
which gives $\{e_0\} = \{e\in \BP_{B^{**}}^\m: e\leq \bp^\ba_B(v)\}$, but this contradicts $e_0\nleq \bp^\ba_B(v)$.

Relation \eqref{eqt:e-0-bot} allows us to pick $e_1\in \BP_{B^{**}}^\m$ satisfying
\begin{equation}\label{eqt:e-1}
e_1e_0=0\quad \text{and} \quad e_1\nleq \bz^\ba_B(v).
\end{equation}
As $B$ is separable, Corollary \ref{cor:Ury}(a) produces an element $a\in \RS_{B^+}$ with
\begin{equation}\label{eqt:a}
\bp_B^\ba(a) = e_1\quad  \text{and} \quad \bz_B^\ba(a) = e_0.
\end{equation}

Our second claim is that $a\notin \{v\}^\de$.
Assume, on the contrary, that $a\in \{v\}^\de$.
By Proposition  \ref{prop:dist-one}(c) and \eqref{eqt:atomic}, one can find $f\in \BP_{B^{**}}^\m$ such that
\begin{equation*}\label{eqt:omega-a-v}
\text{either} \quad f\leq \bp_B^\ba(a) \wedge \bz_B^\ba(v) \quad  \text{or} \quad f\leq \bz_B^\ba(a) \wedge \bp_B^\ba(v).
\end{equation*}
If $f\leq \bp_B^\ba(a) \wedge \bz_B^\ba(v)$, then $f = e_1$ (as $\bp_B^\ba(a) = e_1\in \BP_{B^{**}}^\m$), but \eqref{eqt:e-1} will produce the contradiction that $f\nleq \bz_B^\ba(v)$.
Therefore, we are left with the possibility that $f\leq \bz_B^\ba(a) \wedge \bp_B^\ba(v)$, and hence $f = e_0$ (as $\bz_B^\ba(a) = e_0\in \BP_{B^{**}}^\m$).
However, in this case, \eqref{eqt:e-0} will give us the contradiction that $f \nleq \bp^\ba_B(v)$.
Consequently, our second claim is verified.

On the other hand, as $e_0\leq \bp^\ba_B(u)$ (see \eqref{eqt:e-0}) and $\bz_B^\ba(a) = e_0$ (see \eqref{eqt:a}), Proposition  \ref{prop:dist-one}(c) implies that $a\in \{u\}^\de$.
Hence, one can obtain $a\in \{u\}^\de\setminus \{v\}^\de$, contradicting the assumption.

These establish that if $\{u\}^\de \subseteq \{v\}^\de$,  then  $\bp_A(u) \leq \bp_A(v)$ and $\bz_A(u) \leq \bz_A(v)$.
\end{proof}

As in the proof above, if $\{u\}^\de \not\subseteq \{v\}^\de$, then one can find a positive norm one element $w$ in the unital $C^*$-subalgebra generated by $u$ and $v$ satisfying $\|w-u\| =1$ and $\|w-v\| < 1$.

\smnoind

The following result gives an algebraic property that is determined by the diametral relation (note that as in the proof of Lemma \ref{lem:bij-isom} below, there is a way to identify $\RS_{A^+}^{-1}$ through the diametral relation).

\begin{cor}\label{cor:v-in-hered}
	Suppose that $A$ is a unital $C^*$-algebra.
	For $u,v\in \RS_{A^+}$, one has $\overline{vAv}\subseteq \overline{uAu}$ (or equivalently, $v\in  \overline{uAu}$) if and only if $\{u\}^\de \cap \RS_{A^+}^{-1}\subseteq \{v\}^\de \cap \RS_{A^+}^{-1}$.
\end{cor}
\begin{proof}
By Remark \ref{rem:q-closed}, the inclusion $\overline{vAv}\subseteq \overline{uAu}$ is the same as $v\in  \overline{uAu}$.
Furthermore, Remark \ref{rem:q-closed} and Lemma \ref{lem:peak-zero}(b) also tell  us that
$$\overline{uAu} = \bs_A(u)\cdd A^{**}\cdd \bs_A(u)\cap A = (1-\bz_A(u))A^{**}(1-\bz_A(u))\cap A.$$
Thus, $v\in \overline{uAu}$ if and only if $v(1-\bz_A(u)) = v$; i.e.
$v\bz_A(u)=0$, which   is then equivalent to $\bz_A(u)\leq \bz_A(v)$, because of Lemma \ref{lem:peak-zero}(d).
Now, Proposition \ref{prop:order-1}(a) ensures that $v\in \overline{uAu}$ if and only if $\{u\}^\de \cap \RS_{A^+}^{-1}\subseteq \{v\}^\de \cap \RS_{A^+}^{-1}$.
\end{proof}

The following result gives another algebraic property that is determined by the diametral relation.
This result generalizes
\cite[Theorem 2.3]{peralta2018characterizing}, \cite[Lemma 2.3(a)]{LN}, and \cite[Lemma 4.1]{LNW}.
However, during the final stage of preparing this manuscript, we were informed by A.M. Peralta that he and his coauthor independently obtained the
 same result, in the more general setting of unital $JB^*$-algebras  (see \cite[Theorem 2.1]{peralta2026}).
Nevertheless, the proof given in \cite{peralta2026} is different from ours.

\begin{cor}\label{cor:met-char-proj}
Let $A$ be a unital $C^*$-algebra, and $a\in \RS_{A^+}$.
Then $a$ is a projection if and only if $(\{a\}^\de)^\de = \{a\}$.
\end{cor}
\begin{proof}
Suppose that $a$ is a projection.
Consider $w\in (\{a\}^\de)^\de$.
By \eqref{eqt:111}, one has
$\{a\}^\de = ((\{a\}^\de)^\de)^\de \subseteq \{w\}^\de.$
It then follows from  Proposition \ref{prop:order-1}(b) that
$$\bp_A(a)\leq \bp_A(w)\quad \text{and} \quad \bz_A(a)\leq \bz_A(w).$$
Since $a$ is a projection, $\bp_A(a) =a$ and $\bz_A(a) = 1-a$.
Thus, Lemma \ref{lem:peak-zero}(f) produces
$$\bp_A(w) = a \quad \text{and} \quad \bz_A(w) = 1-a.$$
Hence, Lemma \ref{lem:peak-zero}(g) tells us that $w$ is a projection, and so, $w = \bp_A(w) = a$.
The converse follows from Corollary \ref{cor:Peral-2.2}.
\end{proof}

\section{Diametral relations determine Jordan $^*$-structures}\label{S:weakT}

In this section, we will give an answer to Problem \ref{prob:weak-Tingley}, concerning bijections that preserve the diametral relations (see \eqref{eqt:pre-diam}), in the case when the $C^*$-algebras are unital and have no two dimensional irreducible $^*$-representation.
Observe that the assumption of not having two dimensional irreducible $^*$-representation is quite mild; e.g.,
a $C^*$-algebra of the form $\BM_n(A)$ for some $n\geq 3$ will have no two dimensional non-degenerate $^*$-representation.

\smnoind

Observe that a bijection $\Lambda:\RS_{A^+}\to \RS_{B^+}$ preserves  the diametral relations, in the sense of \eqref{eqt:pre-diam}, if and only if
\begin{equation}\label{eqt:pre-de}
\Lambda(E)^\de = \Lambda(E^\de)\qquad (E\subseteq \RS_{A^+})
\end{equation}
(see \eqref{eqt:def-orth-comp}).
Equality \eqref{eqt:pre-de} will be used implicitly in the following.

\smnoind

Notice also that a bijection that preserves the diametral relations may not be continuous.
The following is such an example.

\begin{eg}\label{eg:counter-eg-lin-ext}
Suppose that both $A$ and $B$ are the two dimensional $C^*$-algebra $\BC\oplus_\infty \BC$.
Define $\Lambda_1:\RS_{A^+}\to \RS_{B^+}$ as
$$\Lambda_1(s,t):= (s^2, 1-t^{1/3}) \qquad (s,t\in [0,1]).$$
Then $\Lambda_1$ is a continuous bijection that preserves the diametral relations, but cannot be extended to a linear map from $A$ to $B$.

On the other hand, the map $\Lambda_2:\RS_{A^+}\to \RS_{B^+}$ defined as
$$
\Lambda_2(s,t):= \begin{cases}
(1-s^{1/2}, t) & \text{when } s,t\in (0,1) \text{ with } t \text{ being rational}\\
(1-s^{1/2}, 1-t) & \text{when } s,t\in (0,1) \text{ with } t \text{ being irrational}\\
(s,t) & \text{when } s,t\in \{0,1\}
\end{cases}
$$
is discontinuous bijection preserving the diametral relations.
\end{eg}

\begin{lem}\label{lem:bij-isom}
	Let $A$ and $B$ be unital $C^*$-algebras.
	Suppose that $\Lambda:\RS_{A^+}\to \RS_{B^+}$ is a bijection preserving the diametral relations (i.e., \eqref{eqt:pre-de} is satisfied).
	
	\smnoind
	(a) $\Lambda(\BP_A\setminus \{0\}) = \BP_B\setminus \{0\}$ and  $\Lambda(\BP_{\Z(A)}\setminus \{0\}) = \BP_{\Z(B)}\setminus \{0\}$.

	\smnoind
	(b) $\Lambda(1) = 1$.
	
	\smnoind
	(c) $\Lambda(\RS_{A^+}^{-1}) = \RS_{B^+}^{-1}$.

\smnoind
(d) For each $u\in \RS_{A^+}$, one has $\Lambda\big(\RS_{(\overline{uAu})^+}\big) = \RS_{(\overline{\Lambda(u)B\Lambda(u)})^+}$.

\smnoind
(e) If $a\in \RS_{A^+}$ and $q\in \BP_A\setminus \{0\}$, then $a\leq q$ if and only if $\Lambda(a)\leq \Lambda(q)$.
\end{lem}
\begin{proof}
(a) This follows from Corollary \ref{cor:met-char-proj} and Lemma \ref{lem:sum-peak}(a).

\smnoind
(b) This follows from part (a) above, together with Lemma \ref{lem:sum-peak}(b).

\smnoind
(c) This follows from part (b) above and the well-known fact that
$$\RS_{A^+}^{-1} = \big\{u\in \RS_{A^+}: \|1-u\| <1 \big\} = \RS_{A^+}\setminus \{1\}^\de.$$
	
\smnoind
(d) Consider $a\in \RS_{A^+}$.
It follows from Corollary \ref{cor:v-in-hered} that $a\in \RS_{(\overline{uAu})^+}$ if and only if $\{u\}^\de \cap \RS_{A^+}^{-1}\subseteq \{a\}^\de \cap \RS_{A^+}^{-1}$.
Moreover, by part (c) above, we have
$$\{u\}^\de \cap \RS_{A^+}^{-1}\subseteq \{a\}^\de \cap \RS_{A^+}^{-1}\quad  \text{if and only if} \quad \{\Lambda(u)\}^\de \cap \RS_{B^+}^{-1}\subseteq \{\Lambda(a)\}^\de \cap \RS_{B^+}^{-1}.$$
Again, one may use Corollary \ref{cor:v-in-hered} to verify that the second inclusion above is equivalent to $\Lambda(a)\in \RS_{(\overline{\Lambda(u)B\Lambda(u)})^+}$.

\smnoind
(e) Notice that $a\leq q$ if and only if $a\in \RS_{(qAq)^+}$.
Since $\Lambda(q)\in \BP_B\setminus\{0\}$ (by part (a) above), the conclusion follows from part (d) above.
\end{proof}

Let us denote $\CC^\ba_\bp(A):= \{\bp^\ba_A(a): a\in \RS_{A^+} \} \subseteq \CC^\ba(A)\setminus \{0\}$ (see \eqref{eqt:peak-non-zero}).
Then
\begin{equation}\label{eqt:CC-a-p}
\CC^\ba_\bp(A)= \big\{\bp^\ba_A(b): b\in \RS_{A^+}^{-1} \big\}  = \big\{\bz^\ba_A(d): d\in \RS_{A^+}\setminus \RS_{A^+}^{-1} \big\} \cup \{1\}.
\end{equation}
In fact, the above follows from the facts that
$$(1+a)/2\in \RS_{A^+}^{-1}\quad \text{and} \quad \bp^\ba_A((1+a)/2) = \bp^\ba_A(a)$$
(see Lemma \ref{lem:peak-zero}(e)),
as well as Lemma \ref{lem:peak-zero}(a),  together with the fact that
when $c\in \RB_{A^+}\setminus \{0\}$, one has (see Lemma \ref{lem:peak-zero}(d))
$$c/\|c\|\in \RS_{A^+}\quad  \text{and} \quad \bz^\ba_A(c/\|c\|) = \bz^\ba_A(c).$$

\begin{lem}\label{lem:inf-peak-proj}
	Let $A$ be a unital $C^*$-algebra.
	
	\smnoind
	(a) For $p,q\in \CC^\ba_\bp(A)$ with $p\wedge q \neq 0$, one has $p\wedge q\in \CC^\ba_\bp(A)$.
	
	\smnoind
	(b) For $q\in \CC^\ba(A)\setminus \{0\}$, there exists a decreasing net $\{b_i\}_{i\in \KI}$ in $\RS_{A^+}^{-1}$ such that $q = \bigwedge_{i\in \KI} \bp^\ba_A(b_i)$, and this implies that  $q = \bigwedge \{p\in \CC^\ba_\bp(A): q\leq p \}$.
\end{lem}
\begin{proof}
(a) Let $a,b\in \RS_{A^+}$ with $\bp^\ba_A(a) = p$ and $\bp^\ba_A(b) =q$.
Set $c:=(a+b)/2$.
For any $e\in \BP_{A^{**}}^\m$, choose $\omega_e\in \ps(A)$ with $\ks_{\omega_e} = e$ (see \eqref{eqt:min-pure-st}).

As $p\wedge q\neq 0$, there is $e_0\in \BP_{A^{**}}^\m[p\wedge q]$  (see \eqref{eqt:atomic}).
By  \eqref{Phi-A-min} and Lemma \ref{lem:peak-zero}(e),
$$e_0\leq p\leq \bp_A(a)\leq a\ \text{ and }\ e_0\leq q\leq \bp_A(b)\leq b,\ \text{ which gives }e_0\leq c.$$
Thus,  $\omega_{e_0}(c) = 1$, and
$c\in \RS_{A^+}.$

Consider $e\in \BP_{A^{**}}^\m[p\wedge q]$.
The same argument as above implies that $e\leq c$, and hence $e\leq \bp^\ba_A(c)$ (see Lemma \ref{lem:peak-zero}(e)); i.e.,
$e \in \BP_{A^{**}}^\m[\bp^\ba_A(c)]$.

Consider $e\in \BP_{A^{**}}^\m[\bp^\ba_A(c)]$.
Then $e\leq \bp^\ba_A(c)\leq \bp_A(c) \leq c$.
Hence, $\omega_e(c) = 1$, which ensures that  $\omega_e(a) = \omega_e(b) = 1$.
From these, $e\leq a$ and $e\leq b$, which is the same as $e\leq p$ and $e\leq q$ (by Lemma \ref{lem:peak-zero}(e)); i.e. $e\in \BP_{A^{**}}^\m[p\wedge q]$.

Now, it follows  from $\BP_{A^{**}}^\m[p\wedge q] = \BP_{A^{**}}^\m[\bp^\ba_A(c)]$ and \eqref{eqt:atomic} that $p\wedge q = \bp^\ba_A(c)$.

\smnoind
(b) By Remark \ref{rem:meet-semilat}(a), Relation \eqref{eqt:Phi-A-closed} and the weak$^*$-continuity of $\Phi_A$,
there is a decreasing net $\{a_i\}_{i\in \KI}$ in $\RS_{A^+}$ with $q = w^*\text{-}{\lim}_i\, \Phi_A(a_i)$.
Lemma \ref{lem:peak-zero}(e) tells us that $q\leq \bp^\ba_A(a_i)\leq \Phi_A(a_i)$ ($i\in \KI$).
Hence,  $q = \bigwedge_{i\in \KI} \bp^\ba_A(a_i)$.
If
$$b_i:= (1+a_i)/2\in \RS_{A^+}^{-1} \qquad (i\in \KI),$$
then Lemma \ref{lem:peak-zero}(e) gives $\bp^\ba_A(b_i) = \bp^\ba_A(a_i)$ ($i\in \KI$), and the conclusion follows.
\end{proof}

\begin{prop}\label{prop:pre-diam-rel}
Let $A$ and $B$ be unital $C^*$-algebras.
Suppose that $\Lambda:\RS_{A^+}\to \RS_{B^+}$ is a bijection preserving the diametral relations.
 Then one can find a bijection $\Gamma:\CC(A)\to \CC(B)$ satisfying
\begin{enumerate}[\upshape{(}i\upshape{)}]
	\item $\Gamma$ is an order isomorphism;
	\item $\Gamma(\bz_A(c)) = \bz_B(\Lambda(c))$ $(c\in \RS_{A^+})$;
	\item $\Gamma$ is bi-orthogonality preserving (see \eqref{eqt:def-bi-ortho-pres});
	\item $\Gamma(\bp_A(c)) = \bp_B(\Lambda(c))$ $(c\in \RS_{A^+})$;
	\item $\Gamma(p) = \Lambda(p)$ $(p\in \BP_A\setminus \{0\})$.
\end{enumerate}
\end{prop}
\begin{proof}
We will first construct  a bijection $\Gamma_0:\CC^\ba(A)\to \CC^\ba(B)$
satisfying
\begin{enumerate}
	\item $\Gamma_0$ is an order isomorphism;
	\item $\Gamma_0(\bz^\ba_A(c)) = \bz^\ba_B(\Lambda(c))$ ($c\in \RS_{A^+}$);
	\item $\Gamma_0$ is bi-orthogonality preserving;
	\item $\Gamma_0(\bp^\ba_A(c)) = \bp^\ba_B(\Lambda(c))$ ($c\in \RS_{A^+}$);
	\item $\Gamma_0(\Phi_A(p)) = \Phi_B(\Lambda(p))$\; ($p\in \BP_A\setminus \{0\}$).
\end{enumerate}
Then Lemma \ref{lem:Phi-A-closed}(a) will produce a bijection $\Gamma:\CC(A)\to \CC(B)$ with
$$\Phi_B\circ \Gamma = \Gamma_0\circ \Phi_A,$$
and Lemma \ref{lem:Phi-A-closed}(b) implies that $\Gamma$ satisfies Conditions (i) - (v).

\smnoind

Before the construction, we note that since $\BP_{A^{**}}^\m\subseteq \CC^\ba(A) \subseteq \BP_{A^{**}_\ba}$ (see Lemma \ref{lem:Phi-A-closed}(c)), it follows from \eqref{eqt:atomic} that for $p\in \CC^\ba(A)\setminus \{0\}$,
\begin{equation}\label{eqt:q-closed-sup}
\text{$p$ is the supremum of $\BP_{A^{**}}^\m[p]$ in $\CC^\ba(A)$}.
\end{equation}
Moreover, we recall from Lemma \ref{lem:Phi-A-closed}(f) that for every family
$\{p_\lambda\}_{\lambda\in \Upsilon}$ in $\CC^\ba(A)$,
\begin{equation}\label{eqt:inf-q-cl}
{\bigwedge}_{\lambda\in \Upsilon}\, p_\lambda\in \CC^\ba(A)\text{ and is the infimum of $\{p_\lambda\}_{\lambda\in \Upsilon}$ in }\CC^\ba(A).
\end{equation}

The construction will be done through several steps.

\smnoind
\textbf{Step I:} We first construct an order isomorphism $\ti \Gamma_0: \CC^\ba_\bp(A)\to \CC^\ba_\bp(B)$.

More precisely, for $p\in \CC^\ba_\bp(A)$ and $a\in \RS_{A^+}^{-1}$ with
$\bp^\ba_A(a) = p$ (see \eqref{eqt:CC-a-p}), we set
\begin{equation}\label{eqt:def-ti-Gamma-0}
\ti \Gamma_0(p) := \bp^\ba_B(\Lambda(a))\in \CC^\ba_\bp(B).
\end{equation}
In order to show that $\ti \Gamma_0$ is a well-defined order preserving map, we
consider $q\in \CC^\ba_\bp(A)$ with $p\leq q$ and $b\in \RS_{A^+}^{-1}$ satisfying $\bp^\ba_A(b)  = q$.
Since $\bp^\ba_A(a)\leq \bp^\ba_A(b)$ and $\bz^\ba_A(a) = \bz^\ba_A(b) = 0$ (see Proposition \ref{prop:dist-one}(d)),
we know from Proposition \ref{prop:order-1}(b) that $\{a\}^\de\subseteq \{b\}^\de$.
It then follows from (see \eqref{eqt:pre-de})
$$
\{\Lambda(a)\}^\de\subseteq \{\Lambda(b)\}^\de
$$
and Proposition \ref{prop:order-1}(b)
ensures that
$$\bp^\ba_B(\Lambda(a)) \leq \bp^\ba_B(\Lambda(b)).$$
It follows that $\ti \Gamma_0$ is both well-defined and order preserving.

By symmetry, there is an order preserving map $\ti \Gamma_0^\#: \CC^\ba_\bp(B)\to \CC^\ba_\bp(A)$ with
\begin{equation*}
\ti \Gamma_0^\#(r) = \bp^\ba_A(\Lambda^{-1}(c)) \quad \text{when} \quad c\in \RS_{B^+}^{-1} \ \ \text{and} \ \ r =\bp_B^\ba(c).
\end{equation*}
If $u\in \RS_{B^+}^{-1}$, then $\Lambda^{-1}(u)\in \RS_{A^+}^{-1}$ (by Lemma \ref{lem:bij-isom}(c)), and
$$\ti \Gamma_0\big(\ti \Gamma_0^\#\big(\bp^\ba_B(u)\big)\big) = \ti \Gamma_0\big(\bp^\ba_A(\Lambda^{-1}(u))\big) = \bp^\ba_B(u).$$
Similarly, $\ti \Gamma_0^\#\circ \ti \Gamma_0 =\id_{\CC^\ba_\bp(A)}$.
This means that $\ti \Gamma_0^\#$ is the inverse of $\ti \Gamma_0$.
Therefore,
\begin{equation}\label{eqt:def-Gamma-1}
	\ti \Gamma_0^{-1}(\bp^\ba_B(u)) = \bp^\ba_A(\Lambda^{-1}(u)) \qquad (u\in \RS_{B^+}^{-1})
\end{equation}
and $\ti \Gamma_0$ is an order isomorphism.

\smnoind
\textbf{Step II:} We will extend $\ti \Gamma_0$ to an order isomorphism $\Gamma_0: \CC^\ba(A) \to \CC^\ba(B)$ such that
\begin{equation}\label{eqt:Theta-min}
	\Gamma_0|_{\BP_{A^{**}}^\m} \text{ is a bijection from }\BP_{A^{**}}^\m \text{ onto }\BP_{B^{**}}^\m.
\end{equation}

For this, we define $\Gamma_0:\CC^\ba(A)\to \CC^\ba(B)$ as $\Gamma_0(0):=0$ and (see \eqref{eqt:inf-q-cl})
\begin{equation}\label{eqt:def-Gamma-0}
\Gamma_0(q) := \bigwedge \{\ti \Gamma_0(p): p\in \CC^\ba_\bp(A); q\leq p \} \qquad (q\in \CC^\ba(A)\setminus \{0\}).
\end{equation}
Since $\ti \Gamma_0$ is order preserving, $\Gamma_0$ is an order preserving map extending $\ti \Gamma_0$.

We claim that
\begin{equation}\label{eqt:Gamma-0-non-zero}
\Gamma_0(q)\neq 0 \qquad (q\in \CC^\ba(A)\setminus \{0\}).
\end{equation}
In fact, if $p_1, \dots, p_n\in \CC^\ba_\bp(A)$ with $q\leq p_k$ for $k=1,\dots,n$, then Lemma \ref{lem:inf-peak-proj}(a) tells us that $p_1\wedge \cdots \wedge p_n\in \CC^\ba_\bp(A)\setminus \{0\}$, and one has $\ti \Gamma_0(p_1)\wedge \cdots \wedge \ti \Gamma_0(p_n) = \ti \Gamma_0(p_1\wedge \cdots \wedge p_n)\neq 0$ (see Step I above and \eqref{eqt:inf-q-cl}).
By Corollary \ref{cor:fin-int-prop}, we see that
$$\Gamma_0(q) = \bigwedge \{\ti \Gamma_0(p): p\in \CC^\ba_\bp(A); q\leq p \} \neq 0.$$

In the same way, if we define a map $\Gamma_0^\#:\CC^\ba(B)\to \CC^\ba(A)$
by $\Gamma_0^\#(0) := 0$ and
\begin{equation}\label{eqt:def-Gamma-0-inv}
 \Gamma_0^\#(s) := \bigwedge \{\ti \Gamma_0^{-1}(r): r\in \CC^\ba_\bp(B); s\leq r \} \qquad (s\in \CC^\ba(B)\setminus \{0\}),
\end{equation}
then $\Gamma_0^\#$ is an order preserving map extending $\ti \Gamma_0^{-1}$ and it satisfies
\begin{equation}\label{eqt:Gamma-0-inv-non-zero}
\Gamma_0^\#(\CC^\ba(B)\setminus \{0\}) \subseteq \CC^\ba(A)\setminus \{0\}.
\end{equation}

Fix $q\in \CC^\ba(A)\setminus \{0\}$.
If $p\in \CC^\ba_\bp(A)$ with $q\leq p$, then \eqref{eqt:def-Gamma-0} implies $\Gamma_0(q)\leq \ti \Gamma_0(p)$.
This means that
$$\ti \Gamma_0 \big(\{p\in \CC^\ba_\bp(A): q\leq p \}\big) \subseteq \{r\in \CC^\ba_\bp(B): \Gamma_0(q)\leq r \}$$
or equivalently (as $\ti \Gamma_0$ is bijective),
$$\{p\in \CC^\ba_\bp(A): q\leq p \} \subseteq \{\ti \Gamma_0^{-1}(r): r\in \CC^\ba_\bp(B);  \Gamma_0(q)\leq r \}.$$
From this, Lemma \ref{lem:inf-peak-proj}(b), \eqref{eqt:Gamma-0-non-zero} and \eqref{eqt:def-Gamma-0-inv}, we have
\begin{equation}\label{eqt:q-Gamma-inv-Gamma-q}
q = \bigwedge \big\{p\in \CC^\ba_\bp(A): q\leq p \big\} \geq \Gamma_0^\#(\Gamma_0(q)).
\end{equation}
Consider $e\in \BP_{A^{**}}^\m$.
We know from \eqref{eqt:q-Gamma-inv-Gamma-q}, \eqref{eqt:Gamma-0-non-zero} and \eqref{eqt:Gamma-0-inv-non-zero} that
\begin{equation*}\label{eqt:Gamma-inv-Gamma-e}
e = \Gamma_0^\#(\Gamma_0(e)) \qquad (e\in \BP_{A^{**}}^\m).
\end{equation*}
Moreover, since both $\Gamma_0$ and $\Gamma_0^\#$ are order-preserving, we know that if $e\leq q$, then
$e = \Gamma_0^\#(\Gamma_0(e))\leq \Gamma_0^\#(\Gamma_0(q))$.
In other words,
$$\{e\in \BP_{A^{**}}^\m: e\leq q\} \subseteq \{e\in \BP_{A^{**}}^\m: e\leq \Gamma_0^\#(\Gamma_0(q))\}$$
and it follows from \eqref{eqt:atomic} that $q\leq \Gamma_0^\#(\Gamma_0(q))$.
This and \eqref{eqt:q-Gamma-inv-Gamma-q} ensure that
$$q = \Gamma_0^\#(\Gamma_0(q)).$$

Similarly, we have $s = \Gamma_0(\Gamma_0^\#(s))$ for every $s\in \CC^\ba(B)$.
These show that $\Gamma_0$ is an order isomorphism with its inverse being $\Gamma_0^\#$.

Equality \eqref{eqt:Theta-min} follows from the fact that minimal projections in $\BP_{A^{**}}$ and $\BP_{B^{**}}$ are minimal elements in $\CC^\ba(A)\setminus \{0\}$ and $\CC^\ba(B)\setminus \{0\}$, respectively.

\smnoind
\textbf{Step III:} We will verify Condition (2).

Let $c\in \RS_{A^+}$.
If $\bz^\ba_A(c) = 0$, then $c\in \RS_{A^+}^{-1}$ (see Proposition \ref{prop:dist-one}(d)) and Lemma \ref{lem:bij-isom}(c), together with Proposition \ref{prop:dist-one}(d), implies that $\bz^\ba_B(\Lambda(c)) = 0 = \Gamma_0(\bz^\ba_A(c))$.

Suppose that $\bz^\ba_A(c) \neq 0$.
Pick any (see \eqref{eqt:atomic})
$$e\in \BP_{A^{**}}^\m[\bz^\ba_A(c)].$$
By Lemma \ref{lem:inf-peak-proj}(b), there exists a decreasing net $\{b_i\}_{i\in \KI}$ of elements in $\RS_{A^+}^{-1}$ satisfying
\begin{equation}\label{eqt:e=inf}
e = {\bigwedge}_{i\in \KI}\, \bp^\ba_A(b_i).
\end{equation}
Proposition \ref{prop:dist-one}(c) ensures that $\|b_i-c\| =1$ for every $i\in \KI$.
Hence,
$$\|\Lambda(b_i)-\Lambda(c)\| =1 \qquad (i\in \KI).$$
For $i\in \KI$, since $\bz^\ba_B(\Lambda(b_i)) = 0$ (by Lemma \ref{lem:bij-isom}(c) and Proposition \ref{prop:dist-one}(d)), we learn from Proposition \ref{prop:dist-one}(c) that $\bp^\ba_B(\Lambda(b_i))\wedge \bz^\ba_B(\Lambda(c))\neq 0$; i.e., (see Step II and \eqref{eqt:def-ti-Gamma-0})
\begin{equation}\label{eqt:Gamma-0-wedge}
\Gamma_0\big(\bp^\ba_A(b_i)\big)\wedge \bz^\ba_B(\Lambda(c)) = \bp^\ba_B(\Lambda(b_i))\wedge \bz^\ba_B(\Lambda(c))\neq 0.
\end{equation}
Pick $i_1, \dots,i_n\in \KI$.
There exists $i_0\in \KI$ with $b_{i_k}\geq b_{i_0}$ for $k=1,\dots,n$ (because $\{b_i\}_{i\in \KI}$ is decreasing).
As $\Gamma_0$ is order preserving, \eqref{eqt:inf-q-cl} and \eqref{eqt:Gamma-0-wedge} produce
$$\Gamma_0\big(\bp^\ba_B(b_{i_1})\big)\wedge  \cdots \wedge \Gamma_0\big(\bp^\ba_B(b_{i_n})\big)\wedge \bz^\ba_B(\Lambda(c)) \geq \Gamma_0\big(\bp^\ba_B(b_{i_0})\big)\wedge \bz^\ba_B(\Lambda(c))\neq 0.$$
Moreover, since $\Gamma_0$ is an order isomorphism, \eqref{eqt:inf-q-cl}, \eqref{eqt:e=inf} and Corollary \ref{cor:fin-int-prop} imply
$$\Gamma_0(e)\wedge \bz^\ba_B(\Lambda(c)) = {\bigwedge}_{i\in \KI} \Gamma_0(\bp^\ba_B(b_i))\wedge \bz^\ba_B(\Lambda(c))\neq 0.$$
Hence, $\Gamma_0(e)\leq \bz^\ba_B(\Lambda(c))$, since $\Gamma_0(e)\in \BP_{B^{**}}^\m$ (by \eqref{eqt:Theta-min}).
This establishes that
$$\Gamma_0\left(\BP_{A^{**}}^\m[\bz^\ba_A(c)] \right) \subseteq \BP_{B^{**}}^\m[\bz^\ba_B(\Lambda(c)) ].$$
By considering $\Lambda^{-1}$, we see from the above argument and \eqref{eqt:def-Gamma-1} (instead of \eqref{eqt:def-ti-Gamma-0}) that
$\Gamma_0^{-1}\left(\BP_{B^{**}}^\m[\bz^\ba_B(\Lambda(c))] \right) \subseteq \BP_{A^{**}}^\m[\bz^\ba_A(c)]$.
These give
$$\Gamma_0\left(\BP_{A^{**}}^\m[\bz^\ba_A(c)]\right) = \BP_{B^{**}}^\m[\bz^\ba_B(\Lambda(c))].$$
Now, because $\Gamma_0: \CC^\ba(A)\to \CC^\ba(B)$ is an order isomorphism, we know from \eqref{eqt:q-closed-sup} that $\Gamma_0$ satisfies Condition (2).

\smnoind
\textbf{Step IV:} We claim that for $e_0,e_1\in \BP_{A^{**}}^\m$, one has
\begin{equation}\label{eqt:orth}
e_0e_1 =0 \quad \text{if and only if} \quad \Gamma_0(e_0)\Gamma_0(e_1) = 0.
\end{equation}

In fact, suppose that $e_0e_1 =0$.
Proposition \ref{prop:non-comm-Ury+fin-int-prop}(a), as well as parts (d) and (e) of Lemma \ref{lem:peak-zero}, produce $a\in \RS_{A^+}$ with
$$e_1\leq \bp_A^\ba(a) \quad \text{and} \quad e_0\leq \bz_A^\ba(a).$$
Set $b:=(1+a)/2\in \RS_{A^+}^{-1}$.
Then $\bp^\ba_A(b) = \bp^\ba_A(a)$ (see Lemma \ref{lem:peak-zero}(e)) and $\bz^\ba_A(b) = 0$ (see Proposition \ref{prop:dist-one}(d)).
From these and Proposition \ref{prop:order-1}(b), one has
$$\{b\}^\de\subseteq \{a\}^\de.$$
Hence,
$\{\Lambda(b)\}^\de\subseteq \{\Lambda(a)\}^\de$
(see \eqref{eqt:pre-de}),
and Proposition \ref{prop:order-1}(b) then implies
$$
\bp^\ba_B(\Lambda(b))\leq \bp^\ba_B(\Lambda(a)).
$$
By the definition of $\Gamma_0$ and the fact that it is order preserving, we know that
$$\Gamma_0(e_1) \leq \Gamma_0(\bp^\ba_A(a)) = \Gamma_0(\bp^\ba_A(b)) = \bp^\ba_B(\Lambda(b)) \leq \bp^\ba_B(\Lambda(a)).$$
Furthermore, since $\Gamma_0$ is order preserving, we learn from Step III that
$$\Gamma_0(e_0)  \leq \Gamma_0(\bz^\ba_A(a)) = \bz^\ba_B(\Lambda(a)).$$
Now, it follows from Lemma \ref{lem:peak-zero}(f) that $\Gamma_0(e_1)\Gamma_0(e_0) = 0$.
Conversely, by considering the inverse $\Gamma_0^{-1}$ of $\Gamma_0$, we know that $\Gamma_0(e_1)\Gamma_0(e_0) = 0$ will imply $e_1e_0 = 0$.
Thus, \eqref{eqt:orth} is valid.

\smnoind
\textbf{Step V:} We will verify Condition (3).

Actually, one can conclude from \eqref{eqt:q-closed-sup}, Step II and \eqref{eqt:orth} that for  $p,q\in \CC^\ba(A)$,
\begin{equation*}\label{eqt:orth2}
	\Gamma_0(p)\Gamma_0(q) = 0 \quad \text{if and only if} \quad pq =0.
\end{equation*}

\smnoind
\textbf{Step VI:} Since $\Gamma_0$ is a bi-orthogonality preserving order isomorphism, we see that
$\Gamma_0^{-1}$ is also a bi-orthogonality preserving order isomorphism.

\smnoind
\textbf{Step VII:}
We claim that
\begin{equation*}
\text{ for $q\in \CC^\ba_\bp(B)$ and $v\in \RS_{B^+}$ with $q\leq \bp^\ba_B(v)$, one has $\Gamma_0^{-1}(q) \leq \bp^\ba_A(\Lambda^{-1}(v))$}
\end{equation*}
(this does not follows directly from \eqref{eqt:def-Gamma-1} because $v$ may not belong to $\RS_{B^+}^{-1}$).

Indeed, \eqref{eqt:CC-a-p} gives $u\in \RS_{B^+}^{-1}$ with $\bp^\ba_B(u)  =q$.
Since $\bz^\ba_B(u) = 0$ (see Proposition \ref{prop:dist-one}(d)) and $\bp^\ba_B(u)\leq \bp^\ba_B(v)$, we learn from Proposition \ref{prop:order-1}(b) that $\{u\}^\de\subseteq \{v\}^\de$.
Hence, $\{\Lambda^{-1}(u)\}^\de \subseteq \{\Lambda^{-1}(v)\}^\de$ (see \eqref{eqt:pre-de}),
which produces (via Proposition \ref{prop:order-1}(b))
$$\bp^\ba_A(\Lambda^{-1}(u)) \leq \bp^\ba_A(\Lambda^{-1}(v)).$$
These, together with \eqref{eqt:def-Gamma-1}, verify the required relation.

\smnoind
\textbf{Step VIII:}
In the same way, we also have
\begin{equation*}
	\text{ for $p\in \CC^\ba_\bp(A)$ and $c\in \RS_{A^+}$ with $p\leq \bp^\ba_A(c)$, one has $\Gamma_0(p) \leq \bp^\ba_B(\Lambda(c))$.}
\end{equation*}

\smnoind
\textbf{Step IX:} We will verify Condition (4) (note that Condition (4) does not follows from the definition of $\Gamma_0$ because $c\in \RS_{A^+}$ may not be invertible).

Assume, on the contrary, that there exists $c\in \RS_{A^+}$ with
$$\Gamma_0(p)\neq \bp^\ba_B(\Lambda(c)), \quad \text{when}\quad p:=\bp^\ba_A(c). $$
By Step VIII, we have $\Gamma_0(p)\leq \bp^\ba_B(\Lambda(c))$.
Since $\Gamma_0(p)\neq \bp^\ba_B(\Lambda(c))$, Remark \ref{rem:proj-atomic}(iii) produces $f\in \BP_{B^{**}_\ba}^\m$ that satisfies
$$f\leq \bp^\ba_B(\Lambda(c))\quad  \text{but}\quad f\Gamma_0(p) = 0.$$
It follows from Step VI and the equality $f\Gamma_0(p) = 0$ that
\begin{equation}\label{eqt:Gamma-1-f-p}
\Gamma_0^{-1}(f) p =0.
\end{equation}
On the other hand, as $f\leq \bp^\ba_B(\Lambda(c))$ and $\Gamma_0^{-1}$ is an order isomorphism (see Step VI), we see from Step VII (with $q= \bp_B^\ba(\Lambda(c))\in \CC^\ba_\bp(B)$ and $v = \Lambda(c)$) that
$$\Gamma_0^{-1}(f)\leq \Gamma_0^{-1}(q) \leq \bp^\ba_A(c) = p,$$
which contradicts \eqref{eqt:Gamma-1-f-p}, because $\Gamma_0^{-1}(f)\neq 0$ (see \eqref{eqt:Theta-min}).
These give Condition (4).

\smnoind
\textbf{Step X:} We will verify Condition (5).

In fact, this follows from Step IX and Lemma \ref{lem:bij-isom}(a).
More precisely, let $p\in \BP_A\setminus \{0\}$.
Then $\bp^\ba_A(p) = \Phi_A(\bp_A(p)) = \Phi_A(p)$.
By Step IX and Lemma \ref{lem:bij-isom}(a),
\begin{align*}
\Gamma_0(\Phi_A(p)) & = \Gamma_0(\bp^\ba_A(p)) = \bp^\ba_B(\Lambda(p)) = \Phi_B(\bp_B(\Lambda(p))) = \Phi_B(\Lambda(p)),
\end{align*}
as required.
\end{proof}

The existence of $\Gamma$ as in the above gives us some connections between the $C^*$-algebras $A$ and $B$.
In particular, we have the following.

\begin{prop}\label{prop:pre-diam-rel-bc}
Let $A$ and $B$ be unital $C^*$-algebras,
and $\Lambda:\RS_{A^+}\to \RS_{B^+}$ be a bijection preserving the diametral relations.
Let $\Gamma:\CC(A)\to \CC(B)$ be the bijection as obtained in Proposition \ref{prop:pre-diam-rel}.
If $C\subseteq A$ is a hereditary $C^*$-subalgebra and $r\in \CC(A)$ with $C = (1-r)A^{**}(1-r)\cap A$ (see \eqref{eqt:bij-open-proj-here}), then we put
$$C^\Lambda :=(1-\Gamma(r))B^{**}(1-\Gamma(r))\cap B.$$

\smnoind
(a) If $C$ is an ideal, then so is $C^\Lambda$.

\smnoind
(b)  $C^\Lambda = \overline{\Lambda(u)B\Lambda(u)}$ when $C = \overline{uAu}$.

\smnoind
(c) $C_1^\Lambda C_2^\Lambda = \{0\}$ when $C_1,C_2\subseteq A$ are hereditary $C^*$-subalgebras  with $C_1C_2 = \{0\}$.

\smnoind
(d) If $u_1,u_2\in \RS_{A^+}$ with $u_1u_2 = 0$, then $\Lambda(u_1)\Lambda(u_2) = 0$.
\end{prop}
\begin{proof}
Let us denote $A[s]:= sA^{**}s \cap A$ when $s\in \BP_{A^{**}}$.

\smnoind
(a) Note that $C$ is an ideal if and only if $r$ is central (see Remark \ref{rem:meet-semilat}(d)).
Therefore, this part follows from Lemma \ref{lem:cent-cl-proj}.

\smnoind
(b) By Remark \ref{rem:q-closed} and Lemma \ref{lem:peak-zero}(b),
$$C = \overline{uAu} = A[\bs_A(u)] = A[1-\bz_A(u)].$$
Hence, Condition (ii) of Proposition \ref{prop:pre-diam-rel}, Lemma \ref{lem:peak-zero}(b) and Remark \ref{rem:q-closed} imply
\begin{align*}
C^\Lambda & = B[1-\Gamma(\bz_A(u))] = B[1-\bz_B(\Lambda(u)]
= B[\bs_B(\Lambda(u)] = \overline{\Lambda(u)B\Lambda(u)}.
\end{align*}

\smnoind
(c) For $i=1,2$, if $r_i\in \CC(A)$ satisfying $C_i = A[1-r_i]$, then
$$C_i^\Lambda = B[1-\Gamma(r_i)],$$
and $\overline{C_i}^{\sigma(A^{**},A^*)} = (1-r_i)A^{**}(1-r_i)$  (see \eqref{eqt:A[r]-bidual}).
Since $C_1C_2=\{0\}$, one has $\overline{C_1}^{\sigma(A^{**},A^*)}\overline{C_2}^{\sigma(A^{**},A^*)}=\{0\}$, which implies $(1-r_1)(1-r_2) = 0$.
If $p_i:= \Phi_A(r_i)$ ($i=1,2$), then $(\ida{A}-p_1)(\ida{A}-p_2) =  0$.
Consequently,
\begin{equation}\label{eqt:min-ortho}
\text{for minimal projections $e_1,e_2$ with $e_1\bot p_1$ and $e_2\bot p_2$, one has $e_1\bot e_2$.}
\end{equation}
Let $\Gamma_0:\CC^\ba(A)\to \CC^\ba(B)$ be the bi-orthogonality preserving order isomorphism as in the proof of Proposition \ref{prop:pre-diam-rel}.
Then $\Gamma_0(\BP_{A^{**}}^\m) = \BP_{B^{**}}^\m$.
Moreover, by \eqref{eqt:min-ortho}, for $f_1,f_2\in \BP_{B^{**}}^\m$ with $f_1\bot \Gamma_0(p_1)$ and $f_2\bot \Gamma_0(p_2)$, one has $f_1\bot f_2$.
Thus, \eqref{eqt:atomic} implies
$$(\ida{B}-\Gamma_0(p_1))(\ida{B}-\Gamma_0(p_2)) = 0.$$
From this, we conclude that $\Phi_B\big(C_1^\Lambda\big)\Phi_B\big(C_2^\Lambda\big) = \{0\}$, and hence $C_1^\Lambda C_2^\Lambda = \{0\}$, because  $\Phi_B|_B$ is injective.

\smnoind
(d) Set $C_1 := \overline{u_1 A u_1}$ and $C_2 := \overline{u_2 A u_2}$.
As $C_1C_2 = \{0\}$, parts (b) and (c) give
$$\overline{\Lambda(u_1)B\Lambda(u_1)}\, \overline{\Lambda(u_2)B\Lambda(u_2)} = C_1^\Lambda C_2^\Lambda = \{0\}.$$
From this, we know that $\Lambda(u_1)\Lambda(u_2) = 0$, since $\Lambda(u_i)\in  \overline{\Lambda(u_i)B \Lambda(u_i)}$.
\end{proof}

It follows from part (b) above and Lemma \ref{lem:bij-isom}(d) that if $A$ is separable, then for every non-zero hereditary $C^*$-subalgebra $C\subseteq A$, one has
$$\Lambda(\RS_{C^+}) = \RS_{(C^\Lambda)^+}.$$
We will see in Proposition \ref{prop:met-pre-here} that  the separable assumption can be dropped if $\Lambda$ is assumed to be continuous.

\smnoind

Lemma \ref{lem:bij-isom} and Proposition \ref{prop:pre-diam-rel} also help us to verify the following facts.
Notice that part (b) below is different from Lemma \ref{lem:bij-isom}(e).

\begin{cor}\label{cor:1-a-b}
	Let $A$ and $B$ be unital $C^*$-algebras.
	Suppose that $\Lambda:\RS_{A^+}\to \RS_{B^+}$ is a bijection preserving the diametral relations.

\smnoind
(a) For every $b\in \RS_{A^+}\setminus \RS_{A^+}^{-1}$ and $c\in \RS_{A^+}$, one has $\|1-b - c\| = 1$ if and only if $\|1 - \Lambda(b) - \Lambda(c)\| = 1$.

\smnoind
(b) For $p\in \BP_A\setminus \{0,1\}$ and $a\in \RS_{A^+}^{-1}$, one has
$p\leq a$ if and only if $\Lambda(p)\leq \Lambda(a)$.

\smnoind
(c) $\Lambda$ restricts to a bi-orthogonality preserving order isomorphism from  $\BP_A\setminus \{0\}$ onto $\BP_B\setminus \{0\}$.

\smnoind
(d) For $p\in \BP_A\setminus \{0,1\}$, one has $\Lambda(1-p) = 1 -\Lambda(p)$.
\end{cor}
\begin{proof}
(a) As $b\notin \RS_{A^+}^{-1}$, we have $1-b\in \RS_{A^+}$.
Moreover, by Lemma \ref{lem:bij-isom}(c), we know that $\Lambda(b)\notin \RS_{B^+}^{-1}$, and $1-\Lambda(b)\in \RS_{B^+}$.
	Note that $\|(1-b) - c\| = 1$ if and only if
	$$\text{$\bz_A(b) \wedge \bz_A(c) \neq 0$ \quad or \quad $\bp_A(b) \wedge \bp_A(c) \neq 0$}$$
	(because of Proposition \ref{prop:dist-one}(c) and Lemma \ref{lem:peak-zero}(a)), which is equivalent to
	$$\text{
		$\bz_B(\Lambda(b)) \wedge \bz_B(\Lambda(c)) \neq 0$  \quad or \quad $\bp_B(\Lambda(b)) \wedge \bp_B(\Lambda(c)) \neq 0$},$$
	(by considering the map $\Gamma$ given by Proposition \ref{prop:pre-diam-rel} and use Remark \ref{rem:inf-closed=closed}(b)), and the latter is equivalent to
	$$\|(1 - \Lambda(b)) - \Lambda(c)\| = 1$$
	(again, by Proposition \ref{prop:dist-one}(c)  and Lemma \ref{lem:peak-zero}(a)).

\smnoind
(b) We know from parts (b) and (c) of Lemma \ref{lem:bij-isom} that $\Lambda(1) =1$ and $\Lambda(\RS_{A^+}^{-1}) = \RS_{B^+}^{-1}$.
On the other hand, Lemma \ref{lem:p-leq-a} tells us that $p\leq a$ if and only if
$$\|1- b - a\| = 1\text{ for every
	$b\in \RS_{A^+}$ with } b\leq p,$$
which, by part (a) above, is the same as (note that $b\notin \RS_{A^+}^{-1}$ because $b\leq p\neq 1$)
$$\|1- \Lambda(b) - \Lambda(a)\| = 1\text{ for every
	$b\in \RS_{A^+}$ with } b\leq p,$$
and this can be rewritten as (because of Lemma \ref{lem:bij-isom}(e))
$$\|1- c - \Lambda(a)\| = 1\text{ for every
	$c\in \RS_{B^+}$ with } c\leq \Lambda(p)$$
which is then equivalent to $\Lambda(p)\leq \Lambda(a)$ (by Lemma \ref{lem:p-leq-a} and Lemma \ref{lem:bij-isom}(a)).

\smnoind
(c) This follows from Lemma \ref{lem:bij-isom}(a) and Proposition \ref{prop:pre-diam-rel}.

\smnoind
(d) By part (c) above, if we extend $\Lambda|_{\BP_A\setminus \{0\}}$ to a bijection $\tilde \Lambda: \BP_A \to \BP_B$ that sends $0$ to $0$, then $\tilde \Lambda$ is bi-orthogonality preserving.
Hence, $\tilde \Lambda(1-p) = 1-\tilde \Lambda(p)$ for every $p\in \BP_A$ (see e.g. \cite[Lemma 1]{dye1955geometry}).
\end{proof}

If, in addition, $\Lambda$ is continuous, then we have a better version than part (b):
\begin{quotation}
for $p\in \BP_A\setminus \{0,1\}$ and $a\in \RS_{A^+}$ with $p\leq a$, one has $\Lambda(p)\leq \Lambda(a)$.
\end{quotation}
The reason is that one can find a decreasing sequence $\{a_n\}_{n\in \BN}$ in $\RS_{A^+}^{-1}$ that norm converge to $a$.

\smnoind

Our next result below gives a solution of Problem \ref{prob:weak-Tingley}, in the case when $A$ has no two dimensional irreducible $^*$-representation.

\smnoind

\begin{thm}\label{thm:pre-diam-rel}
	Let $A$ and $B$ be unital $C^*$-algebras such that $A$ has no $2$-dimensional irreducible $^*$-representation.
	If there exists a bijection $\Lambda:\RS_{A^+}\to \RS_{B^+}$ preserving the diametral relations,
then $A$ and $B$ are Jordan $^*$-isomorphic under a Jordan $^*$-isomorphism $\Psi:A\to B$ satisfying
	\begin{itemize}\itemsep=6pt
		\item $\Psi^{**}(\bp_A(a)) = \bp_B(\Lambda(a)) \qquad (a\in \RS_{A^+})$;	
		\item $\Psi^{**}(\bz_A(a)) = \bz_B(\Lambda(a)) \qquad (a\in \RS_{A^+})$;
		\item $\Psi|_{\BP_A\setminus \{0\}} = \Lambda|_{\BP_A\setminus \{0\}}$.
	\end{itemize}
\end{thm}
\begin{proof}
Let $\Gamma:\CC(A)\to \CC(B)$ be the bijection as in Proposition \ref{prop:pre-diam-rel}.
Since $\Gamma$ is bi-orthogonality preserving,
Proposition \ref{prop:DN-3.8} produces a $^*$-isomorphism $\Psi:A\to B$ such that $\Gamma = \Psi^{**}|_{\CC(A)}$.
From Conditions (ii), (iv) and (v) of Proposition \ref{prop:pre-diam-rel}, one knows that $\Psi^{**}$ satisfies the three conditions in the statement.
\end{proof}

Theorem \ref{thm:pre-diam-rel} tells us that  if $A$ has no $2$-dimensional irreducible $^*$-representation, then for every bijection $\Lambda:\RS_{A^+}\to \RS_{B^+}$ that preserves the diametral relations, its restriction $\Lambda|_{\BP_A\setminus \{0\}}$ extends to a Jordan $^*$-isomorphism from $A$ to $B$.
Nevertheless, as we have seen in Example \ref{eg:counter-eg-lin-ext}, it may not be possible to extend $\Lambda$ to a linear map.

\smnoind

In the next section, we will show that if $\Lambda$ is metric preserving, then $\Lambda$ itself extends to a Jordan $^*$-isomorphism from $A$ to $B$, in the case when $A$ has real rank zero (which could have $2$-dimensional irreducible $^*$-representation).

\section{Bijective isometries between positive unit spheres}\label{S:T}

The aim of this section is to give a positive answer to Problem \ref{prob:Tingley} in the case when the algebras are unital and have real rank zero.

\smnoind

In the following, we set
$$A[a]:= \overline{aAa} \quad \text{and} \quad A[r]:= rA^{**}r \cap A \qquad (a\in \RB_{A^+}, r\in \BP_{A^{**}});$$
i.e, $A[a]$ is the smallest hereditary $C^*$-subalgebra of $A$ containing $a$.
Note that when $p\in \BP_A$, the two definitions of $A[p]$ above coincide, and they agree with $pAp$.

\begin{prop}\label{prop:met-pre-here}
	Let $A$ and $B$ be unital $C^*$-algebras, and $\Lambda:\RS_{A^+} \to \RS_{B^+}$ be a continuous bijection preserving the diametral relations.
If $C\subsetneq A$ is a non-zero hereditary $C^*$-subalgebra and $C^\Lambda\subseteq B$ is the hereditary $C^*$-subalgebra as in Proposition \ref{prop:pre-diam-rel-bc}, then $\Lambda(\RS_{C^+}) = \RS_{(C^\Lambda)^+}$.
\end{prop}
\begin{proof}
Let $\{a_i\}_{i\in \KI}$ be an approximate unit of $C$.
Denote $C_i:=A [a_i]$.
Since $\{a_i\}_{i\in \KI}$ is increasing, $\bigcup_{i\in \KI} C_i$ is a norm-dense $^*$-subalgebra of $C$.
Lemma \ref{lem:bij-isom}(d) and Proposition \ref{prop:pre-diam-rel-bc}(b) give
\begin{equation}\label{eqt:Lambda-RS-Ci+}
\Lambda(\RS_{C_i^+}) = \RS_{(C_i^\Lambda)^+}.
\end{equation}
Consider $r\in \CC(A)\setminus \{0,1\}$ with $C = A[1-r]$ (see \eqref{eqt:bij-open-proj-here}).
Then
\begin{equation}\label{eqt:C-Lambda}
C^\Lambda =  B[1- \Gamma(r)],
\end{equation}
where $\Gamma:\CC(A)\to \CC(B)$ is the order isomorphism as in Proposition \ref{prop:pre-diam-rel}.
As $\{a_i\}_{i\in \KI}$ is an increasing net that weak$^*$-converges to $1-r$ (see \eqref{eqt:A[r]-bidual}), we learn from
$$1-r \geq \bs_A(a_i) \geq a_i \qquad (i\in \KI)$$
(by Lemma \ref{lem:peak-zero}(c)) that
$1-r = \bigvee_{i\in \KI} \bs_A(a_i)$, and  Lemma \ref{lem:peak-zero}(b) implies that
$$r = {\bigwedge}_{i\in \KI} \bz_A(a_i).$$
This implies that $r$ is the infimum of $\{\bz_A(a_i)\}_{i\in \KI}$ in $\CC(A)$.
Since $\Gamma$ is an order isomorphism, $\Gamma(r)$ is the infimum of $\{\Gamma(\bz_A(a_i))\}_{i\in \KI}$ in $\CC(B)$.
Consequently,
\begin{equation*}
\Gamma(r) = {\bigwedge}_{i\in \KI} \Gamma(\bz_A(a_i))
\end{equation*}
(see Remark \ref{rem:inf-closed=closed}(b)).
Since $\{\bz_A(a_i)\}_{i\in \KI}$ is decreasing (by Lemma \ref{lem:peak-zero}(d)),
\begin{equation}\label{eqt:Gamma-p=inf}
\text{the increasing net $\{1-\Gamma(\bz_A(a_i))\}_{i\in \KI}$ will weak$^*$-converge to }1- \Gamma(r)
\end{equation}
On the other hand, Remark \ref{rem:q-closed} and Lemma \ref{lem:peak-zero}(b) imply
$$C_i = A[\bs_A(a_i)]= A[1-\bz_A(a_i)] \qquad (i\in \KI).$$
From this, we obtain
$$C_i^\Lambda = B[1-\Gamma(\bz_A(a_i))] \qquad (i\in \KI).$$
This, together with \eqref{eqt:C-Lambda} and \eqref{eqt:Gamma-p=inf}, ensures that $\bigcup_{i\in \KI} C_i^\Lambda$ is a weak$^*$-dense, and hence norm-dense, $^*$-subalgebra of $C^\Lambda$.
Now, the conclusion follows from \eqref{eqt:Lambda-RS-Ci+} and a continuity argument.
\end{proof}
\smnoind

We have seen in Example \ref{eg:counter-eg-lin-ext}, there exists continuous bijection that preserves the diametral relations, but does not preserve the metric.
From now on, we will consider metric preserving bijections between positive unit spheres.

\begin{lem}\label{lem:met-pre-proj-central}
	Let $A$ and $B$ be unital $C^*$-algebras, and $\Lambda:\RS_{A^+} \to \RS_{B^+}$ be a surjective isometry.

	\smnoind
	(a) For each $p\in \BP_A\setminus \{0,1\}$, if we define $\bar \Delta_{[1-p]}: \RB_{{A[1-p]}^+} \to B$ by
	\begin{align}\label{eqt:def-bar-Delta}
		\bar \Delta_{[1-p]}(a) :=	\Lambda(p  + a) - \Lambda(p) \quad \big(a\in \RB_{{A[1-p]}^+}\big),
	\end{align}
	then $\bar \Delta_{[1-p]}$ extends to  Jordan $^*$-isomorphism $\Delta_{[1-p]}: A[1-p] \to B[1-\Lambda(p)]$
	with
	\begin{equation}\label{eqt:Lambda0a-2}
	\Delta_{[1-p]}(r) = \Lambda(r)\qquad (r\in \BP_{A[1-p ]}\setminus \{0\}).
	\end{equation}
	
	\smnoind
	(b) If $p_0,p_1,\dots,p_m\in \BP_A\setminus \{0,1\}$ are orthogonal and $\alpha_1,\dots, \alpha_m\in (0,1]$, then
	\begin{equation*}\label{eqt:Lambda-conv-proj}
		\Lambda\Big(p_0 + {\sum}_{k=1}^{m} \alpha_k p_k\Big) = \Lambda(p_0) + {\sum}_{k=1}^{m} \alpha_k \Lambda(p_k).
	\end{equation*}
%
\end{lem}
\begin{proof}
	(a) Let us set $q:= \Lambda(p)$.
	Corollary \ref{cor:1-a-b}(b) implies that for each $a\in \RS_{A^+}^{-1}$,
	\begin{equation*}\label{eqt:proj-dom-inver}
		\text{$p\leq a$}  \quad \text{if and only if}\quad q\leq \Lambda(a).
	\end{equation*}	
	From this, together with Lemma \ref{lem:bij-isom}(c) and  Lemma \ref{lem:peak-and-zero}(d), one has
	$$
	\Lambda\big(p + \RB_{{A[1-p]}^+}^{-1}\big) = q  + \RB_{{B[1-q ]}^+}^{-1}.
	$$
	As $\Lambda$ is an isometry, we learn from the norm density of $\RB_{{A[1-p ]}^+}^{-1}$ in $\RB_{{A[1-p ]}^+}$ that
	$$
	\Lambda\big(p +\RB_{{{A[1-p ]}}^+}\big) = q  +\RB_{{B[1-q ]}^+},
	$$
	and this gives a surjective isometry $\bar \Delta_{[1-p]}: \RB_{{A[1-p ]}^+} \to \RB_{{B[1-q ]}^+}$ satisfying \eqref{eqt:def-bar-Delta}. 	
	Since $\bar \Delta_{[1-p]}(0) = 0$ and $\RB_{{A[1-p ]}^+}$ is the translation of the closed ball of $A[1-p]_\sa$ with radius $1/2$ and centre $0$, it follows from \cite[Theorem 5]{mankiewicz1972extension} that $\bar \Delta_{[1-p]}$ extends to a real linear bijective isometry $\Delta_{[1-p]}: A[1-p ]_\sa \to B[1-q ]_\sa$ sending $1-p$ to $1-q$.
	The complexification of $\Delta_{[1-p]}$, again denoted by $\Delta_{[1-p]}$,  is then a Jordan $^*$-isomorphism from $A[1-p]$ onto $B[1-q]$ (see, e.g.,  \cite{Kad52}).

	To verify \eqref{eqt:Lambda0a-2}, we note that by \eqref{eqt:def-bar-Delta},
	\begin{equation}\label{eqt:Lambda-Phi}
	\Lambda(p +r) = \Lambda(p ) + \Delta_{[1-p]}(r).
	\end{equation}
	Since $\Lambda|_{\BP_A\setminus \{0\}}$ is an orthogonality preserving order isomorphism from $\BP_A\setminus \{0\}$ onto $\BP_B\setminus \{0\}$ (see Corollary \ref{cor:1-a-b}(c)) and $pr=0$, we know that
	$$\Lambda(p +r) = \Lambda(p \vee r)  = \Lambda(p )\vee \Lambda(r) = \Lambda(p )+\Lambda(r).$$
	We conclude from this and \eqref{eqt:Lambda-Phi} that $\Lambda(r) = \Delta_{[1-p]}(r)$.

\smnoind
(b) It follows from \eqref{eqt:def-bar-Delta} and \eqref{eqt:Lambda0a-2} that
	\begin{align*}\label{eqt:Lambda-Phi}
		\Lambda\Big(p_0 + {\sum}_{k=1}^m \alpha_k p_k\Big)  & =  \Lambda(p_0) + \Delta_{[1-p_0]}\Big({\sum}_{k=1}^m \alpha_k p_k\Big) \\
		& = \Lambda(p_0) + {\sum}_{k=1}^m \alpha_k \Delta_{[1-p_0]}(p_k) = \Lambda(p_0) +   {\sum}_{k=1}^m \alpha_k \Lambda(p_k),
	\end{align*}
as required.
\end{proof}

\begin{lem}\label{lem:ext}
Let $A$ and $B$ be unital $C^*$-algebras such that $A$ has real rank zero.
If $\Lambda:\RS_{A^+} \to \RS_{B^+}$ is a surjective isometry such that $\Lambda|_{\BP_A\setminus \{0\}}:\BP_A\setminus \{0\} \to \BP_B\setminus \{0\}$ extends to a bounded linear map $\Theta: A \to B$, then $\Theta|_{\RS_{A^+}} = \Lambda$ and $\Theta(A) = B$.
\end{lem}
\begin{proof}
By Lemma \ref{lem:met-pre-proj-central}(b), we know that
$$\Lambda\Big(p_0 + {\sum}_{k=1}^{n} \alpha_k p_k\Big) = \Theta\Big(p_0 + {\sum}_{k=1}^{n} \alpha_k p_k\Big),$$
for pairwise orthogonal elements  $p_0,p_1,\dots,p_n\in \BP_A\setminus \{0,1\}$ and $\alpha_1,\dots, \alpha_n\in (0,1]$.
Since $A$ has real rank zero, any element $a\in \RS_{A^+}$ can be approximate by elements of the form $p_0 + {\sum}_{k=1}^{n} \alpha_k p_k$ as above.
The continuity of both $\Lambda$ and $\Theta$ tell us that
$$\Lambda(a) = \Theta(a)\qquad (a\in \RS_{A^+}).$$
We then know from $\Theta(\RS_{A^+}) = \Lambda(\RS_{A^+}) = \RS_{B^+}$ that $\Theta(A) = B$.
\end{proof}

\begin{prop}\label{prop:no-irred-2-dim}
	Let $A$ and $B$ be unital $C^*$-algebras such that $A$ has real rank zero but has no $2$-dimensional irreducible $^*$-representation.
	If $\Lambda:\RS_{A^+} \to \RS_{B^+}$ is a surjective isometry, then $\Lambda$ extends to a Jordan $^*$-isomorphism from $A$ onto $B$.
\end{prop}
\begin{proof}
  By Theorem \ref{thm:pre-diam-rel}, the bijection  $\Lambda|_{\BP_A\setminus \{0\}}:\BP_A\setminus \{0\} \to \BP_B\setminus \{0\}$ extends to a Jordan $^*$-isomorphism $\Theta: A \to B$.
	The conclusion follows from Lemma \ref{lem:ext}.
\end{proof}


The above verifies Theorem \ref{thm:main2}, in the case when $A$ ``has no type $\I_2$ part''.
It is well-known that any von Neumann algebra is a direct sum of a type $\I_2$ algebra and an algebra without type $\I_2$ part.
However, for a general $C^*$-algebra $A$, one can only find an ideal $I$ without $2$-dimensional irreducible $^*$-representation such that $A/I$ is a $2$-homogenous $C^*$-algebra.
Due to the difficulty in dealing with extension, we will not prove Theorem \ref{thm:main2} via such a ``decomposition''.
To present our proof for Theorem \ref{thm:main2}, we need the following result.

\begin{lem}\label{lem:decomp-proj}
	Let $A$ and $B$ be unital $C^*$-algebras such that $A$ has
 real rank zero. Suppose that
 $\Lambda:\RS_{A^+} \to \RS_{B^+}$ is a surjective isometry.
For each $r\in \BP_A\setminus \{0,1\}$, we denote by $\Delta_{[r]}$ the Jordan $^*$-isomorphism as in Lemma \ref{lem:met-pre-proj-central}(a).

\smnoind
(a)	If $r\in \BP_A\setminus \{0,1\}$, then
	$\Lambda|_{\RS_{A_r^+}} = (\Delta_{[r]}\oplus \Delta_{[1-r]})|_{\RS_{A_r^+}}$, where
$$
A_r:= A[r] + A[1-r].
$$	

\smnoind
(b)	If $C\subsetneq A$ is a non-trivial hereditary $C^*$-subalgebras and $C^\Lambda$ is as in Proposition \ref{prop:met-pre-here}, then there is a Jordan $^*$-isomorphism $\Delta_{C}: C \to C^\Lambda$ with  $\Delta_{C}|_{\RS_{C^+}} = \Lambda|_{\RS_{C^+}}$.
\end{lem}
\begin{proof}
Note that $\Delta_{[r]}$ is a Jordan $^*$-isomorphism from $rAr$ to $\Lambda(r)B\Lambda(r)$, because of Corollary \ref{cor:1-a-b}(d).
Moreover, Corollary \ref{lem:bij-isom}(d) gives
$$\Lambda(\RS_{(rAr)^+}) = \RS_{(\Lambda(r)B\Lambda(r))^+}.$$

\smnoind
(a) Let $\alpha_1,\dots, \alpha_m, \beta_1,\dots, \beta_n$ be elements in  $(0,1]$ satisfying
\begin{equation}\label{eqt:max-1}
	\max \{ \alpha_1,\dots, \alpha_m, \beta_1,\dots, \beta_n\} = 1.
\end{equation}
and $p_1,\dots, p_m, q_1,\dots, q_n\in \BP_A\setminus \{0,1\}$ are pairwise orthogonal with $p_1,\dots, p_m\leq r$ and $q_1,\dots, q_n\leq 1-r$.
By \eqref{eqt:max-1}, Lemma \ref{lem:met-pre-proj-central}(b) as well as  \eqref{eqt:Lambda0a-2}, one has
\begin{align*}
	\Lambda\Big({\sum}_{k=1}^{m} \alpha_k p_k + {\sum}_{k=1}^{n} \beta_l q_l\Big) & = {\sum}_{k=1}^{m} \alpha_k \Lambda(p_k)
	+ {\sum}_{l=1}^{n} \beta_k \Lambda(q_l)\\
	& = {\sum}_{k=1}^{m} \alpha_k \Delta_{[r]}(p_k)
	+ {\sum}_{l=1}^{n} \beta_k \Delta_{[1-r]}(q_l)\\
	& = \Delta_{[r]}\Big({\sum}_{k=1}^{m} \alpha_k p_k\Big) + \Delta_{[1-r]}\Big({\sum}_{k=1}^{n} \beta_l q_l\Big).
\end{align*}
Since $\Lambda|_{\RS_{A_r^+}}$ and $(\Delta_{[r]} \oplus \Delta_{[1-r]})|_{\RS_{A_r^+}}$ are continuous and elements in $\RS_{A_r^+}$ can be approximated by elements of the form
$${\sum}_{k=1}^m \alpha_k p_k + {\sum}_{l=1}^n \beta_l q_l,$$
as above, we obtain the required equality $\Lambda|_{\RS_{A_r^+}} = (\Delta_{[r]}\oplus \Delta_{[1-r]})|_{\RS_{A_r^+}}$.

\smnoind
(b) It follows from \cite[Corollary 2.8]{BP91} that  $C$ has real rank zero.
We first consider the case when there exists $a\in \RS_{A^+}\setminus \RS_{A^+}^{-1}$ with
$$C = A[a].$$
It is easy to see that $a$ is a strictly positive element of $C$.
Thus, $C$ is $\sigma$-unital, and \cite[Proposition 2.9]{BP91} produces an increasing sequence $\{r_n\}_{n\in \BN}$ in $\BP_C\setminus \{0\}\subseteq \BP_A\setminus \{0,1\}$ that is an approximate unit for $C$.

For $n \in\BN$,  Lemma \ref{lem:ext} and \eqref{eqt:Lambda0a-2} imply that $\Delta_{[r_n]}|_{\RS_{(r_n A r_n)^+}} = \Lambda|_{\RS_{(r_n A r_n)^+}}$.
Hence,
$$
\Delta_{[r_{n+1}]}|_{r_n A r_n} = \Delta_{[r_n]} \qquad (n\in\BN).
$$
If we set $C_0:= {\bigcup}_{n\in \BN} A[r_n]$, then $C_0$ is a norm-dense $^*$-subalgebra of $C$, and
the above produces an isometric Jordan $^*$-homomorphism
$\Delta_0:   C_0 \to B$ with
\begin{equation*}
	\Delta_0|_{r_n A r_n} = \Delta_{[r_n]} \quad \text{and} \quad \Delta_0|_{\RS_{(r_n A r_n)^+}} = \Lambda|_{\RS_{(r_n A r_n)^+}}\qquad (n\in \BN).
\end{equation*}
Thus, $\Delta_0$ extends to an isometric Jordan $^*$-homomorphism
$\Delta_C: C\to B$
such that
\begin{align}\label{eq:DLagree}
\Delta_C|_{\RS_{C^+}} = \Lambda|_{\RS_{C^+}}.
\end{align}
Proposition \ref{prop:met-pre-here} then implies that $\Delta_C(C) = C^\Lambda$.

For the general case, consider $\{a_i\}_{i\in \KI}$ to be an approximate unit of $C$.
Since $C\neq A$, none of the $a_i$ is invertible.
Moreover, $\bigcup_{i\in \KI} A[a_i]$  is a norm-dense $^*$-subalgebra of $C$.
Hence, the first case above and its argument will produce the required Jordan $^*$-isomorphism.
\end{proof}

\begin{prop}\label{prop:non-simple}
Let $A$ and $B$ be unital $C^*$-algebras.
Suppose that $A$ has real rank zero and admits a $2$-dimensional irreducible $^*$-representation.
Then every surjective isometry $\Lambda:\RS_{A^+} \to \RS_{B^+}$ extends to a Jordan $^*$-isomorphism from $A$ onto $B$.
\end{prop}
\begin{proof}
First of all, we consider the case when $A = \BM_2$.
By Proposition \ref{prop:pre-diam-rel}, there is an order isomorphism $\Gamma :\BP_{\BM_2} = \CC(\BM_2)\to \CC(B)$ extending $\Lambda|_{\BP_{\BM_2}\setminus \{0\}}$.
Since
$$\Gamma(\CC(\BM_2)\setminus \{0\}) = \Lambda(\BP_{\BM_2}\setminus \{0\}) = \BP_B\setminus \{0\}$$ (see Lemma \ref{lem:bij-isom}(a)), we have
$\CC(B) = \BP_B.$
Thus, \eqref{eqt:min>closed} implies that $\BP_{B^{**}}^\m \subseteq \BP_B$.
From this and Corollary \ref{cor:1-a-b}(c), for each $p\in \BP_{B^{**}}^\m$, there exists at most one $q\in \BP_{B^{**}}^\m$ that is orthogonal to $p$.
Since rank-one projections from different irreducible $^*$-representations are orthogonal, we see that $B$ either have only one irreducible $^*$-representation and it is $2$-dimensional, or have only two irreducible $^*$-representations, each of them is $1$-dimensional.
Since $\BP_B\setminus \{0,1\}$ is infinite (because $\BP_{\BM_2}\setminus \{0,1\}$ is infinite), we conclude that $B\cong \BM_2$.
Now, one may apply Winger's theorem  (see, e.g., \cite[Theorem 1.1]{BJM}) to extend $\Lambda$ to a Jordan $^*$-isomorphism from $A$ onto $B$.

It remains to consider the case when there exists a non-injective irreducible $^*$-repre\-sentation $\pi:A\to \BM_2$.

Denote  $I_\pi:=\ker \pi\neq \{0\}$.
It follows from \eqref{eqt:bij-open-proj-here} and Remark \ref{rem:meet-semilat}(d) that there is a unique central closed projection $z\in \CC(A)\cap \Z(A^{**})\setminus \{0\}$ satisfying
\begin{equation*}
	I_\pi :=  A[ 1-z].
\end{equation*}
If $z\in \BP_A$, then $A\cong (1-z)A \oplus \BM_2$ and the conclusion follows from Lemma \ref{lem:decomp-proj}(a).
Therefore, we assume that
\begin{equation}\label{eqt:z-not-in-A}
	z\notin \BP_A.
\end{equation}
We will establish the required conclusion through several steps.

In the following, for a norm closed subspace $E\subseteq A$, we identify $E^{**}$ with the weak$^*$-closure $\overline{E}^{\sigma(A^{**}, A^*)}\subseteq A^{**}$ in the canonical way.

\smnoind
\textbf{Step I:}
For a non-trivial hereditary $C^*$-subalgebra $C\subsetneq A$, we denote by
$$\Delta_C:C\to C^\Lambda$$
 the Jordan $^*$-isomorphism given by Lemma \ref{lem:decomp-proj}(b), which satisfies
\begin{equation}\label{eqt:Delta=Lambda}
\Delta_{C}|_{\RS_{C^+}} = \Lambda|_{\RS_{C^+}}.
\end{equation}
Let
$\Delta_C^{**}: C^{**} \to (C^\Lambda)^{**}$
be the normal Jordan $^*$-isomorphism extending $\Delta_C$.
Suppose that $D\subsetneq A$ is another hereditary $C^*$-subalgebra with $D\cap C\neq \{0\}$.
Then as $\Delta_{C}|_{\RS_{(D\cap C)^+}} = \Lambda|_{\RS_{(D\cap C)^+}} = \Delta_{D}|_{\RS_{(D\cap C)^+}}$,
we know that
\begin{equation}\label{eqt:intersect}
\Delta_C^{**}|_{(D\cap C)^{**}} = \Delta_{D}^{**}|_{(D\cap C)^{**}}.
\end{equation}

We recall from Proposition \ref{prop:pre-diam-rel-bc} that if $C = A[1-r]$ for some $r\in \CC(A)$, then
$$C^\Lambda = B[1-\Gamma(r)].$$
In particular, $I_\pi^\Lambda = B[1-\Gamma(z)]$ is an ideal (see Proposition \ref{prop:pre-diam-rel-bc}(a)) and
we obtain a normal Jordan $^*$-isomorphism
$$
\Delta_{I_\pi}^{**}: (1-z) A^{**} \to \big(1-\Gamma(z)\big)B^{**}
$$
 (see \eqref{eqt:A[r]-bidual} and Remark \ref{rem:meet-semilat}(d)).

\smnoind
\textbf{Step II:} \cite[Corollary 2.8]{BP91} and \eqref{eqt:z-not-in-A} allow us to pick $q_0\in \BP_{I_\pi}\setminus \{0\}$ with
\begin{equation}\label{eqt:q0-not-1-z}
q_0\lneqq 1-z.
\end{equation}
We learn from $z\leq 1-q_0$ and \eqref{eqt:A[r]-bidual} that
\begin{equation}\label{eqt:z-c}
	zA^{**} \subseteq (1-q_0)A^{**} = A[1-q_0]^{**}.
\end{equation}
As $q_0\neq 0$, we know that $A[1-q_0]\neq A$.
By Condition (v) of Proposition \ref{prop:pre-diam-rel},
\begin{equation}\label{eqt:Gamma-q0}
\Gamma(q_0) = \Lambda(q_0).	
\end{equation}
Consider $\Delta_{A[1-q_0]}:A[1-q_0]\to A[1-\Lambda(q_0)]$ to be the injective Jordan $^*$-homo\-morphism as in Step I above, and we set
\begin{equation}\label{eqt:Psi}
	\bar \Xi_{q_0}:=\Delta_{A[1-q_0]}^{**}|_{zA^{**}}.
\end{equation}

\smnoind
\textbf{Step III:} We define a normal Jordan $^*$-homomorphism $\Theta: A^{**} \to B^{**}$ by
\begin{equation}\label{eqt:def-Theta}
	\Theta := \Delta_{I_\pi}^{**} \oplus \bar \Xi_{q_0}
\end{equation}
(where $q_0$ be as in Step II), and claim that $\Theta(1) = 1$.

Indeed, as in Step I, one has
$$
\Theta(1-z) = \Delta_{I_\pi}^{**}(1-z) = 1 -\Gamma(z).
$$
On the other hand, we learn from $q_0\bot z$, Condition (iii) of Proposition \ref{prop:pre-diam-rel} as well as  \eqref{eqt:Gamma-q0} that $\Gamma(z)\bot \Lambda(q_0)$.
Since $q_0 + z = q_0 \vee z\in \CC(A)$ (by Proposition \ref{prop:non-comm-Ury+fin-int-prop}(b)), Condition (i) of Proposition \ref{prop:pre-diam-rel} and \eqref{eqt:Gamma-q0} imply
$$\Gamma(q_0 +  z) = \Gamma(q_0 \vee  z) = \Lambda(q_0) + \Gamma(z).$$
Hence, Step I produces a Jordan $^*$-isomorphism (see \eqref{eqt:q0-not-1-z})
$$\Delta_{A[1-q_0 - z]}: A[1-q_0 - z]\to B[1-\Lambda(q_0) - \Gamma(z)]$$
such that $\Delta_{A[1-q_0 - z]}^{**} = \Delta_{A[1-q_0]}^{**}|_{A[1-q_0 - z]^{**}}$ (see \eqref{eqt:intersect}).
It then follows from
$$\Delta_{A[1-q_0]}^{**}(1-q_0) = 1-\Lambda(q_0) \ \ \text{and} \ \ \Delta_{A[1-q_0 - z]}^{**}(1-q_0-z) = 1-\Lambda(q_0) - \Gamma(z)$$
that $1 -\Lambda(q_0) - \Delta_{A[1-q_0]}^{**}(z)= \Delta_{A[1-q_0]}^{**}(1-q_0-z) = 1-\Lambda(q_0) - \Gamma(z)$; i.e.
$$\Theta(z) = \Delta_{A[1-q_0]}^{**}(z)= \Gamma(z)$$
as required.

\smnoind
\textbf{Step IV:}
We will verify that
\begin{equation}\label{eqt:Theta-restr-here}
\Theta|_{C^{**}} = \Delta_{C}^{**}, \text{ for each hereditary $C^*$-subalgebra  $C\subsetneq A$ that contains $I_\pi$}.	
\end{equation}

In fact, \eqref{eqt:bij-open-proj-here} produces $s\in \CO(A)\setminus \{1\}$ with
$$C = A[s].$$
Since $I_\pi\subseteq C$, one has $1-z\leq s$.
As $z$ is central, we see from \eqref{eqt:A[r]-bidual} that
$$C^{**} = s A^{**} s = (1-z) A^{**} (1-z) \oplus sz A^{**} sz = I_\pi^{**}\oplus sz A^{**} sz.$$
By \eqref{eqt:def-Theta}  and \eqref{eqt:intersect}, we know that
$$\Theta|_{I_\pi^{**}} = \Delta_{I_\pi}^{**} = \Delta_C^{**}|_{I_\pi^{**}}.$$

It remains to show that
$$\Theta|_{szA^{**} sz} = \Delta_{C}^{**}|_{szA^{**} sz}.$$
Indeed, without loss of generality, assume that $sz\neq 0$.
It follows from $q_0\leq 1-z\leq s$ that $s(1-q_0) = s - q_0 = (1-q_0)s$.
Proposition \ref{prop:non-comm-Ury+fin-int-prop}(b) then ensures that $s\wedge (1-q_0)\in \CO(A)$ and
one has
\begin{equation}\label{eqt:A-r-wedge}
A[s\wedge (1-q_0)]\subseteq A[s]\cap A[1-q_0].
\end{equation}
On the other hand, we see from
$$s-q_0 \geq s - (1-z) = s - s(1-z) = sz$$
that $sz\leq 1-q_0$, which gives $sz\leq s\wedge (1-q_0)$.
Therefore, \eqref{eqt:A[r]-bidual} and \eqref{eqt:A-r-wedge} give
\begin{align*}
szA^{**}sz
& \subseteq (s\wedge (1-q_0))A^{**}(s\wedge (1-q_0))\\
& = A[s\wedge (1-q_0)]^{**} \subseteq (A[s]\cap A[1-q_0])^{**}.
\end{align*}
Thus, we learn from \eqref{eqt:intersect} that
$$\Theta|_{szA^{**} sz} = \bar \Xi_{q_0}|_{szA^{**} sz} =\Delta_{A[1-q_0]}^{**}|_{szA^{**} sz} = \Delta_{A[s]}^{**}|_{szA^{**} sz},$$
as required.

\smnoind
\textbf{Step V:} Let $e\in \BP_{\BM_2}$ be a rank one projection.
Then $A_e:= \pi^{-1}(e \BM_2 e)$
is a non-trivial hereditary $C^*$-subalgebra containing $I_\pi$.
Step IV and \eqref{eqt:Delta=Lambda} tell us that
$$\Theta|_{\BP_{A_e}\setminus \{0\}} = \Delta^{**}_{A_e}|_{\BP_{A_e}\setminus \{0\}} = \Lambda|_{\BP_{A_e}\setminus \{0\}}. $$

\smnoind
\textbf{Step VI:} We claim that $\Theta(p) = \Lambda(p)$ ($p\in \BP_{A}\setminus \{0\}$).

Indeed, if $\pi(p) = 0$, then $p\in I_\pi$, and \eqref{eqt:Delta=Lambda} gives $\Theta(p) = \Delta_{I_\pi}(p) = \Lambda(p)$.
If $p =1$, then Lemma \ref{lem:bij-isom}(b) and Step III above give $\Lambda(p) = 1 = \Theta(p)$.
If $\pi(p) = 1$ but $p\neq 1$, then $\pi(1-p) = 0$, and Corollary \ref{cor:1-a-b}(d) as well as the above implies
$$1 - \Lambda(p) = \Lambda(1-p) = \Theta(1-p) = 1 - \Theta(p),$$
which produces $\Theta(p) = \Lambda(p)$.
If $e:=\pi(p)$ is a rank one projection, then Step V tells us that
$\Theta(p) = \Lambda(p)$.
These verified Step VI.

\smnoind

Finally, we learn from Step VI that $\Theta(A) \subseteq B$.
Thus, Lemma \ref{lem:ext} produces $\Theta(A) = B$ and $\Theta|_{\RS_{A^+}} = \Lambda$, as required.
\end{proof}

One of the difficulty in the proof above comes from the fact that one cannot consider $A[z]$ directly (see, e.g., Example \ref{eg:cl-proj-id}).

\smnoind

Let us summarize Proposition \ref{prop:no-irred-2-dim} and Proposition \ref{prop:non-simple} in part (a) of the following result.
Part (b) of this result comes from part (a) and  \cite[Theorem H]{Hers}.

\begin{thm}\label{thm:main}
Let $A$ and $B$ be  unital $C^*$-algebras such that $A$ has real rank zero.
Suppose that there is a surjective isometry $\Lambda:\RS_{A^+} \to \RS_{B^+}$.

\smnoind
(a) $\Lambda$  extends to a Jordan $^*$-isomorphism from $A$ onto $B$.

\smnoind
(b) If, in addition, $A$ is prime (in particular, simple), then $\Lambda$ extends to either a complex linear $^*$-isomorphism or 
a complex linear $^*$-anti-isomorphism from $A$ onto $B$.
\end{thm}

As a corollary, if the requirement  in Theorem \ref{thm:main} holds,  then the existence of a surjective isometry $\Lambda$ as above will imply that $B$ has real rank zero as well.
Moreover, if $A$ is prime (respectively, simple), then $B$ is prime (respectively, simple).

\appendix\section{open/closed projections and q-open/q-closed projections}\label{sec:open-proj}

In this section, we will recall some basic information concerning open/closed projections as well as q-open/q-closed projections.

\smnoind

Throughout this appendix, $A$ is a unital $C^*$-algebra.
For a subspace $E\subseteq A$, we identify $E^{**}$ with the weak$^*$-closure $\overline{E}^{\sigma(A^{**},A^*)}\subseteq A^{**}$.

\smnoind

Let $\pi_u:A\to \CL(\CH_u)$ be the universal representation  of $A$ (i.e., the direct sum of GNS constructions of all states on $A$).
It is well-known that $\pi_u$ extends to a faithful normal unital $^*$-representation of $A^{**}$ with its image being the strong closure of $\pi_u(A)$.
We will forget $\pi_u$ and its normal extension and identity
\begin{align*}
	A\subseteq A^{**}\subseteq \CL(\CH_u).
\end{align*}

\begin{defn}\label{defn:open}(\cite[\S 3.11.10]{Ped})
	A projection $r\in \BP_{A^{**}}$ is called an \emph{open projection} of $A$ if there is  an increasing net $\{a_i\}_{i\in \KI}$ of elements in $A_+$ that strongly-converges to $r$ in $\CL(\CH_u)$.
	Moreover, $s\in \BP_{A^{**}}$ is called a \emph{closed projection} of $A$ if $1-s$ is an open projection  of $A$.
	We denote by
	$$\CO(A)\quad \text{and} \quad \CC(A)$$
	the set of open projections and the set of closed projections of $A$, respectively.
\end{defn}

\begin{rem}\label{rem:meet-semilat}
	Let $r\in \BP_{A^{**}}$ and $A[r]$ be the hereditary $C^*$-subalgebra $r A^{**}r \cap A$.

\smnoind
(a) The following are equivalent:
\begin{enumerate}[\ \ (O1).]
	\item $r\in \CO(A)$;
	
	\item there is an increasing net $\{u_i\}_{i\in \KI}$ in $A_+$ that weak$^*$-converges to $r$;
	
	\item $r\in \overline{A[r]}^{\sigma(A^{**},A^*)}$.
\end{enumerate}
In fact, the equivalence of (O1) and (O2) above is obvious, while the equivalence of (O1) and (O3) follows from \cite[\S 3.11.10]{Ped} and the fact that the strong closure of $A[r]$ equals its weak$^*$-closure.

	\smnoind
	(b) Let $B\subseteq A$ be a $C^*$-subalgebra.
	Consider an approximate unit $\{b_i\}_{i\in \KI}$ of $B$ that  weak$^*$-converges to an element $s\in A^{**}_+$.
	Then $s\in \BP_{A^{**}}$ and
	$\overline{B}^{\sigma(A^{**},A^*)} \subseteq s\cdd A^{**} \cdd s.$
If $B$ is a hereditary $C^*$-subalgebra, then
	$\overline{B}^{\sigma(A^{**},A^*)} = s\cdd A^{**} \cdd s,$
	and hence,
	\begin{equation}\label{eqt:B=A[s]}
	B = \overline{B}^{\|\cdot\|} = \overline{B}^{\sigma(A^{**},A^*)}\cap A =  A[s].
	\end{equation}


	\smnoind
	(c) The equivalence of (O1) and (O3) implies that $r\in \CO(A)$ if and only if
	\begin{equation}\label{eqt:A[r]-bidual}
	\overline{A[r]}^{\sigma(A^{**},A^*)} = r\cdd A^{**}\cdd r.
	\end{equation}
	Hence, it follows from part  (b) that
	\begin{equation}\label{eqt:bij-open-proj-here}
	\text{$s\mapsto A[s]$ is a bijection from $\CO(A)$ onto hereditary $C^*$-subalgebras of $A$.}
	\end{equation}
Consequently, $r\in \CO(A)$ if and only if
\begin{enumerate}[\ \ (O4).]
	\item there is  a hereditary $C^*$-subalgebra $B\subseteq A$ with
	$\overline{B}^{\sigma(A^{**}, A^*)} = r\cdd A^{**}\cdd r.$
\end{enumerate}

\smnoind
(d) By \eqref{eqt:B=A[s]}, an open projection $r\in \CO(A)$ is central if and only if $A[r]$ is an ideal.
\end{rem}

In contrast to the above, when $s\in \CC(A)$, the  $C^*$-subalgebra $A[s]$ may not tell us much information about $s$.

\begin{eg}\label{eg:cl-proj-id}
	Let $A = C([0,1])$ and $x\in [0,1]$.
	Then $I_x:= \{f\in A: f(x) = 0\}$ is an essential ideal of $A$.
	Consider $r_x$ to be the central open projection with $I_x=A[r_x]$.
	As $A^{**} = A^{**}r_x \oplus A^{**}(1-r_x)$ and $A/I_x\cong \BC$, we see that $A^{**}(1-r_x) \cong \BC$.
However, $A[1-r_x] = A\cap A^{**}(1-r_x) = \{0\}$.
\end{eg}

\begin{rem}\label{rem:inf-closed=closed}
(a) Denote $\mathcal{QS}(A):= \{\omega\in A^*_+: \|\omega\|\leq 1 \}$.
	It is well-known that
	$$p\mapsto F(p):= \{\omega\in \mathcal{QS}(A):\omega(p) =1 \}$$
	is an order isomorphism from $\BP_{A^{**}}$ onto the set of norm-closed faces of $\mathcal{QS}(A)$.
	On the other hand, as in \cite[\S 3.11.10]{Ped}, the restriction of the assignment above
	is a bijection from $\CC(A)$ onto the set of weak$^*$-closed faces of $\mathcal{QS}(A)$.

\smnoind
(b) As intersections of weak$^*$-closed faces are weak$^*$-closed, we see from part (a) above that for any subset $E\subseteq \CC(A)$, one has $\bigwedge E \in \CC(A)$ (see also \cite[Theorem II.5]{Ake69}), and so is the infimum of $E$ in $\CC(A)$.
\end{rem}

Let us also recall the following well-known fact (see, e.g., \cite[Theorem 1.5.2]{Ped}).

\begin{lem}\label{lem-left-id-here}
The assignment $L\mapsto L\cap L^*$ is a bijection from
the set of closed left ideals of $A$ onto the set of hereditary $C^*$-subalgebras of $A$.
\end{lem}

For a state $\omega$ of $A$ with support projection $\ks_\omega\in \BP_{A^{**}}$ (see \eqref{eqt:def-s-omega}), the left ideal
$$
L_\omega= \{a\in A: \omega(a^*a)=0\}
$$
of $A$
will satisfy $\overline{L_\omega}^{\sigma(A^{**}, A^*)} = A^{**}\cdd (1-\ks_\omega)$, and
$$\overline{L_\omega^* L_\omega}^{\sigma(A^{**}, A^*)} = (1-\ks_\omega) \cdd A^{**}\cdd (1-\ks_\omega).$$
This and the equivalence of (O1) and (O4) tell us that $\ks_\omega\in \CC(A)$.

\smnoind

In particular, since minimal projections are support projections of pure states,
\begin{equation}\label{eqt:min>closed}
	\BP_{A^{**}}^\m\subseteq \CC(A)
\end{equation}
(see also \cite[Corollary 2]{Akemann68}).
We recall from \eqref{eqt:def-Phi-A} the truncation map $\Phi_A:A^{**}\to A^{**}_\ba$.

\begin{defn}\label{defn:q-open}(\cite{Ake71})
A projection $p\in \BP_{A^{**}_\ba}$ is called a \emph{q-open projection} of $A$ if
$$
A^{**}_\ba\cdd p = \overline{\Phi_A(L)}^{\sigma(A^{**}_\ba, (A^{**}_\ba)_*)}
$$
for a left ideal $L\subseteq A$.
Furthermore, $q\in \BP_{A^{**}_\ba}$ is called a \emph{q-closed projection} of $A$ if $\ida{A}-q$ (see \eqref{eqt:def-id-atom-bidual}) is a q-open projection  of $A$.
We denote by
$$\CO^\ba(A)\quad \text{and} \quad \CC^\ba(A)$$
the set of q-open projections and the set of q-closed projections of $A$, respectively.
\end{defn}

By Lemma \ref{lem-left-id-here} and the equivalence of (O1) and (O4), a projection $r\in \BP_{A^{**}}$ belongs to $\CO(A)$ if and only if there exists a closed left ideal $L\subseteq A$ such that
$$\overline{L}^{\sigma(A^{**}, A^*)} = A^{**}\cdd r.$$
Since $\overline{\Phi_A(L)}^{\sigma(A^{**}_\ba, (A^{**}_\ba)_*)} = \overline{L}^{\sigma(A^{**}, A^*)}\cdd \ida{A}$ for any left ideal $L\subseteq A$, one has
\begin{equation}\label{eqt:Phi-A-closed}
	\Phi_A(\CO(A)) = \CO^\ba(A), \quad \text{and hence} \quad \Phi_A(\CC(A)) = \CC^\ba(A).
\end{equation}

\begin{lem}\label{lem:Phi-A-closed}
Let $B$ be  the $JB$-subalgebra generated by $\CO(A)$.

\smnoind
	(a) The restriction of $\Phi_A$ to $B$ is an injective Jordan homomorphism.
	
	\smnoind
	(b) For $x,y\in A_\sa \cup \CO(A)\cup \CC(A)$, one has $x,y\in B$, and
	\begin{itemize}
		\item $x\leq y$ if and only if $\Phi_A(x) \leq \Phi_A(y)$;
		
		\item $xy=yx$ if and only if $\Phi_A(x) \Phi_A(y) = \Phi_A(y) \Phi_A(x)$.
		
		\item $xy=0$ if and only if $\Phi_A(x) \Phi_A(y) = 0$.
	\end{itemize}
	
	\smnoind
	(c) $\BP_{A^{**}}^\m\subseteq \CC^\ba(A)$.
	
	\smnoind
	(d) For $r\in \CC(A)$, the projection $\Phi_A(r)$ is central if and only if $r$ is central.

\smnoind
(e) If $\{r_\lambda\}_{\lambda\in \Upsilon}$ is a family in $\BP_{A^{**}}$, then $\bigwedge_{\lambda\in \Upsilon} \Phi_A(r_\lambda) = \Phi_A\big({\bigwedge}_{\lambda\in \Upsilon} r_\lambda\big)$.

\smnoind
(f) If $\{p_\lambda\}_{\lambda\in \Upsilon}$ is a family in $\CC^\ba(A)$, then $\bigwedge_{\lambda\in \Upsilon} p_\lambda\in \CC^\ba(A)$ and it is the infimum of $\{p_\lambda\}_{\lambda\in \Upsilon}$ in $\CC^\ba(A)$.
\end{lem}
\begin{proof}
	(a) Denote $A_\sa^m\subseteq A^{**}_\sa$ to be the set of all strong-limits of increasing nets in $A_\sa$, and $\CU$ to be the set of
all ``universally measurable elements'' in $A^{**}_\sa$ (see \cite[\S 4.3.11]{Ped}).
	By \cite[Theorem 4.3.15]{Ped}, the restriction of $\Phi_A$ to $\CU$ is isometric.
	
	Let $\CU_b\subseteq \CU$ be the subspace defined as in \cite[p.437]{Ped72} that contains $A_\sa^m$.
	Then $\CO(A)\subseteq \CU_b$.
	By \cite[Lemma 3.5]{Ped72}, $\CU_b$ is a Jordan subalgebra of $A^{**}_\sa$.
	On the other hand, \cite[Proposition 4.3.13]{Ped} tells us that $\CU$ is norm closed and hence it contain the norm closure of $\CU_b$.
	Consequently, $\CU$ contains $B$.

	\smnoind
	(b) As $\CO(A)\subseteq B$, one has $\CC(A)\subseteq B$.
	Hence, $B$ contains $p-q$, for $p\in \CO(A)$ and $q\in \CC(A)$.
	An approximation argument then implies that $A_\sa\subseteq B$.
	Now, as $\Phi_A|_B:B\to \Phi_A(B)$ is an isometric Jordan isomorphism, the conclusions of this part follows
	(observe that for $a,b\in A^{**}_\sa$, one has
	$ab = ba$ if and only if the Jordan subalgebra generated by $a$ and $b$ is associative).

	\smnoind
	(c) This follows from the second equality in \eqref{Phi-A-min}, together with \eqref{eqt:min>closed} and \eqref{eqt:Phi-A-closed}.
	
	\smnoind
	(d) Clearly, if $r$ is central, then so is $\Phi_A(r)$.
	The converse follows from part (b) above and the weak$^*$-density of $A_\sa$  in $A^{**}_\sa$.

\smnoind
(e) The follows from the fact that $\Phi_A:A^{**}\to A^{**}_\ba$ is the truncation map.

\smnoind
(f) By \eqref{eqt:Phi-A-closed}, there is a family $\{r_\lambda\}_{\lambda\in \Upsilon}$ in $\CC(A)$ such that $p_\lambda = \Phi_A(r_\lambda)$ ($\lambda\in \Upsilon$).
It then follows from part (e) above,  Remark \ref{rem:inf-closed=closed}(b) and \eqref{eqt:Phi-A-closed} that
$${\bigwedge}_{\lambda\in \Upsilon}\, p_\lambda = \Phi_A\Big({\bigwedge}_{\lambda\in \Upsilon} \, r_\lambda\Big) \in \CC^\ba(A).$$
From this, we see that ${\bigwedge}_{\lambda\in \Upsilon}\, p_\lambda$ is the infimum of $\{p_\lambda\}_{\lambda\in \Upsilon}$ in $\CC^\ba(A)$.
\end{proof}

\begin{lem}\label{lem:cent-cl-proj}
	Let $A$ and $B$ be unital $C^*$-algebras.
	Suppose that $\Gamma:\CC(A)\to \CC(B)$ is an order isomorphism that  is also bi-orthogonality preserving.
Then $\Gamma(r)$ is a central closed projection of $B$, when $r\in \CC(A)$ is central.
\end{lem}
\begin{proof}
	By parts (a) and (b) of Lemma \ref{lem:Phi-A-closed} as well as \eqref{eqt:Phi-A-closed}, the map $\Gamma$ induces a bi-orthogonality preserving order isomorphism $\Gamma_\ba:\CC^\ba(A)\to \CC^\ba(B)$.
	Note that a projection $q\in \BP_{A^{**}_\ba}$ is central if and only if for every $e\in \BP_{A^{**}}^\m$, either $e\leq q$ or $e\bot q$ (by \eqref{eqt:atomic}).
	Moreover, it follows from Lemma \ref{lem:Phi-A-closed}(c) that $\BP_{A^{**}}^\m$ is the set of minimal elements in  $\CC^\ba(A)\setminus \{0\}$.
	As $\Gamma_\ba$ is a bi-orthogonality preserving order isomorphism, we see that $\Gamma_\ba(\BP_{A^{**}}^\m) = \BP_{B^{**}}^\m$.
	These imply that $q\in \CC^\ba(A)$ is central if and only if $\Gamma_\ba(q)$ is central.
	The conclusion then follows from Lemma \ref{lem:Phi-A-closed}(d).
\end{proof}


\begin{prop}\label{prop:non-comm-Ury+fin-int-prop}
(a) For any $r,s\in \CC(A)$ with $rs=0$, there is $a\in \RS_{A^+}$ satisfying $ar = r$ and $as = 0$.

\smnoind
(b) If $r,s\in \CC(A)$ with $rs=sr$, then $r\vee s\in \CC(A)$; or equivalently, for $r,s\in \CO(A)$ with $rs = sr$, one has $r\wedge s\in \CO(A)$.
\end{prop}

In fact, part (a) above is taken from \cite[Theorem I.1]{Ake70}, while part (b) is a special case of \cite[Theorem II.7]{Ake69}.

\begin{cor}\label{cor:fin-int-prop}
If $\{p_\lambda\}_{\lambda\in \Upsilon}$ is a family in $\CC^\ba(A)$ with $\bigwedge_{\lambda\in E} p_\lambda \neq 0$ for each non-empty finite subset $E\subseteq \Upsilon$, then $\bigwedge_{\lambda\in \Upsilon} p_\lambda\neq 0$.
\end{cor}
\begin{proof}
By  \eqref{eqt:Phi-A-closed}, there is a family  $\{r_\lambda\}_{\lambda\in \Upsilon}$ in $\CC(A)$ with $p_\lambda = \Phi_A(r_\lambda)$ ($\lambda\in \Upsilon$).
Parts (a) and (b) of Lemma \ref{lem:Phi-A-closed}, Remark \ref{rem:inf-closed=closed}(b) and the hypothesis imply $\bigwedge_{\lambda\in E} r_\lambda \neq 0$, for any non-empty finite subset $E\subseteq \Upsilon$.
Thus, we learn from
\cite[Proposition II.10]{Ake69} that
${\bigwedge}_{\lambda\in \Upsilon} r_\lambda \neq 0.$
By parts (a) and (e) Lemma \ref{lem:Phi-A-closed}, we know that
$${\bigwedge}_{\lambda\in \Upsilon}\, p_\lambda = \Phi_A\Big({\bigwedge}_{\lambda\in \Upsilon}\,  r_\lambda\Big)\neq 0$$
as claimed.
\end{proof}

\section*{Acknowledgement}

Ng is partially supported by the Nankai Zhide Foundation and by National Natural Science Foundation of China (12671155).
Wong is partially supported by the Taiwan NSTC grant 114-2115-M-110-002-MY2.

The authors would like to thank A.M. Peralta for informing them of that an extension of Corollary \ref{cor:met-char-proj} can be found in \cite{peralta2026}.

No result in this article is obtained with the help of AI, although we used AI to help smoothing the Abstract and the Introduction.

\end{document}